\documentclass[preprint,3p,12pt]{elsarticle}
\usepackage{amsmath}
\usepackage{stmaryrd}
\usepackage{bbding}
\usepackage{dcolumn}
\usepackage{graphicx}
\usepackage{amsfonts}
\usepackage{amssymb}
\usepackage{psfrag}
\usepackage{wrapfig}
\usepackage{subfigure}
\usepackage{makeidx}
\usepackage{bm}
\usepackage{epsf}
\usepackage{epsfig}
\usepackage{setspace}
\usepackage{color}
\usepackage{booktabs}
\usepackage{setspace}
\usepackage[ruled,linesnumbered]{algorithm2e}

\definecolor{SpringGreen4}{RGB}{0,139,69}

\begin{document}
\title{Implicit high-order gas-kinetic schemes for compressible flows on three-dimensional unstructured meshes II: unsteady flows}

\author[BNU,HKUST1]{Yaqing Yang}
\ead{yqyangbnu@163.com}

\author[BNU]{Liang Pan\corref{cor}}
\ead{panliang@bnu.edu.cn}

\author[HKUST1,HKUST2]{Kun Xu}
\ead{makxu@ust.hk}

\cortext[cor]{Corresponding author}
\address[BNU]{Laboratory of Mathematics and Complex Systems, School of Mathematical Sciences, Beijing Normal University, Beijing, China}
\address[HKUST1]{Department of Mathematics, Hong Kong University of Science and Technology, Clear Water Bay, Kowloon, Hong Kong}
\address[HKUST2]{Shenzhen Research Institute, Hong Kong University of Science and Technology, Shenzhen, China}
 
\begin{abstract}
For the simulations of unsteady flow, the global time step becomes really small with a large variation of local cell size.  
In this paper,  an implicit high-order gas-kinetic scheme (HGKS) is developed to remove the restrictions on the time step for unsteady simulations. 
In order to improve the efficiency and keep the high-order accuracy, a two-stage third-order implicit time-accurate discretization is proposed. 
In each stage, an artificial steady solution is obtained for the implicit system with the pseudo-time iteration.
In the iteration, the classical implicit methods are adopted to solve the nonlinear system, including the lower-upper symmetric Gauss-Seidel (LUSGS) and generalized minimum residual  (GMRES) methods.
To achieve the spatial accuracy, the HGKSs with both non-compact and compact reconstructions are constructed. 
For the non-compact scheme, the weighted essentially non-oscillatory (WENO) reconstruction is used.  For the compact one, 
the Hermite WENO (HWENO) reconstruction is adopted due to the updates of both cell-averaged flow variables and their derivatives.
The expected third-order temporal accuracy is achieved with the two-stage temporal discretization. 
For the smooth flow, only a single artificial iteration is needed. For uniform meshes, 
the efficiency of the current implicit method improves significantly in comparison with the explicit one.
For the flow with discontinuities, compared with the well-known Crank-Nicholson method, 
the spurious oscillations in the current schemes are well suppressed. 
The increase of the artificial iteration steps introduces extra reconstructions associating with a reduction of the computational efficiency. 
Overall, the current implicit method leads to an improvement in efficiency over the explicit one in the cases with a large variation of mesh size. 
Meanwhile, for the cases with strong discontinuities on a uniform mesh, the efficiency of the current method is comparable with that of the explicit scheme.
\end{abstract}
 

\maketitle

\section{Introduction}
The numerical simulations of compressible unsteady flow require high-order accurate and efficient numerical methods to resolve a wide range of time and length scales. 
Explicit methods have been mainly adopted in such simulations to advance the compressible Navier-Stokes equations in time, and the multi-stage Runge-Kutta method \cite{TVD-RK}  
is usually used as temporal discretization. For unsteady flow computations,  a global time step is used, which is determined by the Courant-Friedrichs-Lewy (CFL) condition. 
For the high speed viscous flow simulation with a significant variation of local cell size, the overall global time step of explicit scheme becomes really small.  
On the other hand, for the problems with the low frequency, such as wing flutter and periodic pitching of an airplane, the time step is also very restricted.  
Thus, it becomes desirable to develop a fully implicit method to release the constraint on the global CFL condition and preserve the stability.
In recent years, significant progress has been made in developing implicit method for the steady flow problems.   
The implicit method requires to solve a large nonlinear system, which arises from the linearization of fully implicit scheme at each time step.  As an approximate factorization method, the lower-upper symmetric Gauss-Seidel (LUSGS) 
method on structured meshes was originally developed by Jameson and Yoon \cite{LUSGS-2}, and 
it has been successfully extended to unstructured meshes \cite{LUSGS-3,LUSGS-4}.
The most attractive feature of this method is that it does not require any extra memory 
in comparison with the explicit one and is free from any matrix inversion. 
To accelerate the convergence,  the nonlinear generalized minimal residual (GMRES) 
method \cite{GMRES-1,GMRES-2} was developed. The GMRES method has a faster convergence speed, 
but the drawback is that it requires a considerable amount of memory to store the Jacobian matrix, 
which may be prohibitive for large problems. To save the storage, the matrix-free GMRES method has been applied to 
unstructured meshes \cite{GMRES-3}. To improve the convergence,
the preconditioners are also applied to the GMRES method \cite{Implicit-2,Implicit-3,Implicit-4,Implicit-5,Implicit-6}.
For unsteady flow problems, the implicit schemes have been developed using dual-time stepping methods \cite{unsteady-implicit-1,unsteady-implicit-2,unsteady-implicit-3}.
In the computation of unsteady flow, an artificial steady solution is obtained for the implicit system with the pseudo-time iteration. For the artificial steady solution, the pseudo-time residual has to be drived to zero (or at least to truncation error).
A fast and matrix-free implicit method has been adopted to obtain the steady solution to the pseudo-time system \cite{GMRES-4}.

In the past decades, the gas-kinetic scheme (GKS) based on the
Bhatnagar-Gross-Krook (BGK) model \cite{BGK-1,BGK-2} has been developed systematically for the computations from subsonic flows
to supersonic ones \cite{GKS-Xu1,GKS-Xu2}. The gas-kinetic scheme
presents a gas evolution process from the kinetic scale to
hydrodynamic scale, and both inviscid and viscous fluxes can be
calculated in one framework.  With the two-stage temporal discretization, 
which was originally developed for the Lax-Wendroff type flow solvers \cite{GRP-high-1}, 
a reliable framework was provided to construct gas-kinetic scheme with
 fourth-order and even higher-order temporal accuracy \cite{GKS-high-1,GKS-high-2}. 
Based on the spatial and temporal coupled property of GKS solver and HWENO reconstruction \cite{HWENO-1,HWENO-2}, 
 the explicit high-order compact gas-kinetic schemes are also developed \cite{GKS-HWENO-1,GKS-HWENO-2}. 
In the compact scheme, the time-dependent gas distribution function at a cell interface is used to calculate the fluxes for the updates of
cell-averaged flow variables, and to evaluate the cell-averaged gradients of flow variables. 
For the gas-kinetic scheme, several implicit gas-kinetic schemes have also been developed for  
the steady flows \cite{GKS-implicit-1,GKS-implicit-2,GKS-implicit-3} and the unsteady 
flows \cite{GKS-implicit-unsteady-1,GKS-implicit-unsteady-2,GKS-implicit-unsteady-3},
which provide efficient techniques for speeding up the convergence. 
Incorporate with the GMRES method, the implicit non-compact scheme and compact scheme are developed for steady problems. 
To further accelerate the computation, the Jacobi iteration is chosen as preconditioner for both non-compact and compact schemes. 
Numerical examples validate the convergence and robustness of implicit HGKS.

For the unsteady flow with a large variation of local mesh size, the global time step is restricted, which reduces the efficiency of computation. 
To remove such restrictions, an implicit HGKS for compressible unsteady flow is developed on the three-dimensional unstructured 
meshes. To achieve the spatial accuracy, the HGKSs with both non-compact and compact reconstructions are constructed. For non-compact HGKS, 
the third-order WENO reconstruction is used,  where the stencils are selected from the neighboring cells and the neighboring cells of 
neighboring cells.  To improve the resolution and parallelism, the compact HWENO reconstruction is used,  where the stencils only contain 
one level of neighboring cells.  In order to improve the computational efficiency and keep a high-order temporal accuracy for unsteady flows,
a two-stage third-order implicit time-accurate discretization is developed.   In each stage, an artificial steady solution is obtained for the implicit system with the pseudo-time iteration. In the iteration, the classical implicit methods are adopted to solve the nonlinear system, including LUSGS and GMRES methods.
Various three-dimensional numerical experiments, from the subsonic to supersonic flows, 
are presented to validate the accuracy and robustness of the current implicit method. The numerical results indicate that the third-order 
accuracy is achieved with the two-stage temporal discretization. 
For the flow with discontinuities, 
compared with the well-known Crank-Nicholson method,
the spurious oscillations in the current schemes are well suppressed. The increase of the artificial iteration steps introduces extra reconstructions associating with a reduction of the computational efficiency. Overall, the
current implicit method leads to an improvement in efficiency over the explicit one in the
cases with a large variation of mesh size.
Meanwhile, for the cases with strong discontinuities on a uniform mesh, the efficiency of the current method is comparable with that of the
explicit scheme.
In the future, the current implicit method will be extended to the moving-mesh problems, in which the time step is usually restricted by the physical frequency.

This paper is organized as follows. In Section 2, BGK equation and finite volume scheme will be introduced. 
The third-order non-compact and compact reconstructions will be presented in Section 3. 
Section 4 includes the time-implicit gas-kinetic scheme for unsteady flow. 
Numerical examples are included in Section 5. The last section is the conclusion.

\section{BGK equation and finite volume scheme}
The Boltzmann equation expresses the behavior of a many-particle kinetic system in terms of the evolution equation for a single
particle gas distribution function. The BGK equation \cite{BGK-1,BGK-2}  is the simplification of Boltzmann equation, and the 
three-dimensional BGK equation can be written as
\begin{equation}\label{bgk}
f_t+uf_x+vf_y+wf_z=\frac{g-f}{\tau},
\end{equation}
where $\boldsymbol{u}=(u,v,w)$ is the particle velocity, $\tau$ is the collision time and $f$ is the gas distribution function. The
equilibrium state $g$ is given by Maxwellian distribution 
\begin{equation*}
g=\rho(\frac{\lambda}{\pi})^{(N+3)/2}e^{-\lambda[(u-U)^2+(v-V)^2+(w-W)^2+\xi^2]},
\end{equation*}
where $\rho$ is the density, $\boldsymbol{U}=(U,V,W)$ is the macroscopic fluid velocity, and $\lambda$ is the inverse of gas temperature, 
i.e., $\lambda=m/2kT$. In the BGK model, the collision operator involves a simple relaxation from $f$ to the local equilibrium state $g$. 
The variable $\xi$ accounts for the internal degree of freedom, $\xi^2=\xi_1^2+\dots+\xi_N^2$, $N=(5-3\gamma)/(\gamma-1)$ is the 
internal degree of freedom, $\gamma$ is the specific heat ratio and takes $1.4$ in the numerical cases. The collision term satisfies the 
compatibility condition
\begin{equation*}
\int \frac{g-f}{\tau}\psi \text{d}\Xi=0,
\end{equation*}
where $\displaystyle\psi=(1,u,v,w,\frac{1}{2}(u^2+v^2+w^2+\xi^2))^T$ and $\text{d}\Xi=\text{d}u\text{d}v\text{d}w\text{d}\xi_1\dots\text{d}\xi_{N}$.
According to the Chapman-Enskog expansion for BGK equation, the macroscopic governing equations can be derived. In the continuum
region, the BGK equation can be rearranged and the gas distribution function can be expanded as
\begin{align*}
f=g-\tau D_{\boldsymbol{u}}g+\tau D_{\boldsymbol{u}}(\tau D_{\boldsymbol{u}})g-\tau D_{\boldsymbol{u}}[\tau D_{\boldsymbol{u}}(\tau D_{\boldsymbol{u}})g]+...,
\end{align*}
where $D_{\boldsymbol{u}}=\displaystyle\frac{\partial}{\partial t}+\boldsymbol{u}\cdot \nabla$. With the zeroth-order truncation $f=g$,  the Euler 
equations can be obtained. For the first-order truncation
\begin{align*}
f=g-\tau (ug_x+vg_y+wg_z+g_t),
\end{align*}
the Navier-Stokes equations can be obtained \cite{GKS-Xu1,GKS-Xu2}.

Taking moments of Eq.\eqref{bgk} and integrating with respect to space, the semi-discretized finite volume scheme can be expressed as
\begin{align}\label{semi}
\frac{\text{d} Q_i}{\text{d} t}=\mathcal{L}(Q_{i}),
\end{align}
where $Q_i=(\rho, \rho U,\rho V, \rho W, \rho E)$ is the cell averaged conservative value of $\Omega_{i}$,  $\rho$ is the density and $\rho E$ is the total energy density. 
The operator $\mathcal{L}$ is defined as
\begin{equation}\label{finite}
\mathcal{L}(Q_{i})=-\frac{1}{|\Omega_{i}|}\sum_{i_p\in N(i)}F_{i,i_p}(t)S_{i_p}=-\frac{1}{|\Omega_{i}|}\sum_{i_p\in N(i)}\iint_{\Sigma_{i_p}}\boldsymbol{F}(Q,t)\text{d}\sigma,
\end{equation}
where $|\Omega_{i}|$ is the volume of $\Omega_{i}$, $\Sigma_{i_p}$ is the common cell interface of $\Omega_{i}$, $S_{i_p}$ is the area of 
$\Sigma_{i_p}$ and $N(i)$ is the set of index for neighboring cells of $\Omega_{i}$. To achieve the expected order of accuracy, the Gaussian 
quadrature is used for the flux integration
\begin{align*}
\iint_{\Sigma_{i_p}}\boldsymbol{F}(Q,t)\text{d}\sigma=\sum_{G}\omega_{G}F_{G}(t)S_{i_p},
\end{align*}
where $\omega_{G}$ is the quadrature weights. The numerical flux $F_{G}(t)$ at Gaussian quadrature
point $\boldsymbol{x}_{G}$ can be given by taking moments of gas
distribution function
\begin{align}\label{flux-G}
F_{G}(t)=\int\psi \boldsymbol{u}\cdot\boldsymbol{n}_{G} f(\boldsymbol{x}_{G},t,\boldsymbol{u},\xi)\text{d}\Xi,
\end{align}
where $F_{G}(t)=(F^{\rho}_{G},F^{\rho U}_{G},F^{\rho V}_{G},F^{\rho W}_{G}, F^{\rho E}_{G})$ and $\boldsymbol{n}_{G}$ is the local normal direction of cell
interface. With the coordinate transformation, the numerical flux in the global coordinate can be obtained. Based on the integral solution of BGK equation 
Eq.\eqref{bgk}, the gas distribution function $f(\boldsymbol{x}_{G},t,\boldsymbol{u},\xi)$ in the local coordinate can be given by
\begin{equation*}
f(\boldsymbol{x}_{G},t,\boldsymbol{u},\xi)=\frac{1}{\tau}\int_0^t
g(\boldsymbol{x}',t',\boldsymbol{u}, \xi)e^{-(t-t')/\tau}\text{d}t'+e^{-t/\tau}f_0(-\boldsymbol{u}t,\xi),
\end{equation*}
where $\boldsymbol{x}'=\boldsymbol{x}_{G}-\boldsymbol{u}(t-t')$ is the trajectory of particles and $f_0$ is the initial gas distribution
function. With the first order spatial derivatives, the second-order gas distribution function at cell interface can be given by
\begin{align}\label{flux}
f(\boldsymbol{x}_{G},t,\boldsymbol{u},\xi)=&(1-e^{-t/\tau})g_0+
((t+\tau)e^{-t/\tau}-\tau)(\overline{a}_1u+\overline{a}_2v+\overline{a}_3w)g_0\nonumber\\
+&(t-\tau+\tau e^{-t/\tau}){\bar{A}} g_0\nonumber\\
+&e^{-t/\tau}g_r[1-(\tau+t)(a_{1}^{r}u+a_{2}^{r}v+a_{3}^{r}w)-\tau A^r)](1-H(u))\nonumber\\
+&e^{-t/\tau}g_l[1-(\tau+t)(a_{1}^{l}u+a_{2}^{l}v+a_{3}^{l}w)-\tau A^l)]H(u),
\end{align}
where the equilibrium state $g_{0}$ and the corresponding
conservative variables $Q_{0}$ can be determined by the
compatibility condition
\begin{align*}
\int\psi g_{0}\text{d}\Xi=Q_0=\int_{u>0}\psi
g_{l}\text{d}\Xi+\int_{u<0}\psi g_{r}\text{d}\Xi.
\end{align*}
With the reconstruction of macroscopic variables, the coefficients
in Eq.\eqref{flux} can be fully determined by the reconstructed
derivatives and compatibility condition
\begin{equation*}
\begin{aligned} 
\displaystyle 
\langle a_{1}^{k}\rangle=\frac{\partial Q_{k}}{\partial \boldsymbol{n_x}}, 
\langle a_{2}^{k}\rangle=\frac{\partial Q_{k}}{\partial \boldsymbol{n_y}}, 
\langle a_{3}^{k}\rangle&=\frac{\partial Q_{k}}{\partial\boldsymbol{n_z}}, 
\langle
a_{1}^{k}u+a_{2}^{k}v+a_{3}^{k}w+A^{k}\rangle=0,\\ 
\displaystyle
\langle\overline{a}_1\rangle=\frac{\partial Q_{0}}{\partial \boldsymbol{n_x}}, 
\langle\overline{a}_2\rangle=\frac{\partial Q_{0}}{\partial \boldsymbol{n_y}},
\langle\overline{a}_3\rangle&=\frac{\partial Q_{0}}{\partial \boldsymbol{n_z}},
\langle\overline{a}_1u+\overline{a}_2v+\overline{a}_3w+\overline{A}\rangle=0,
\end{aligned}
\end{equation*}
where $k=l$ and $r$,  $\boldsymbol{n_x}$, $\boldsymbol{n_y}$, $\boldsymbol{n_z}$ are the unit directions of local coordinate at $\boldsymbol{x}_{G}$
 and $\langle...\rangle$ are the moments of the
equilibrium $g$ and defined by
\begin{align*}
\langle...\rangle=\int g (...)\psi \text{d}\Xi.
\end{align*}
More details of the gas-kinetic scheme can be found in \cite{GKS-Xu1,GKS-Xu2}.

\section{High-order spatial reconstruction}
To deal with the complex geometry, the three-dimensional tetrahedral and
hexahedral unstructured meshes are considered. In the previous studies, the high-order gas-kinetic schemes have 
been developed with the third-order non-compact WENO reconstruction \cite{GKS-WENO-3} and compact 
HWENO reconstruction \cite{GKS-HWENO-1, GKS-HWENO-2}.  Successes have been achieved for the unsteady flows 
with explicit schemes from the subsonic to supersonic flow problems, and a brief review is given in this section.

\subsection{Non-compact WENO reconstruction}
For the target cell $\Omega_i$, the faces are labeled as $F_{p}$, $p=1,\dots,P$, where $P=4$ for tetrahedral cell and $P=6$ for hexahedral cell. 
The neighboring cell of $\Omega_{i}$, which shares the face $F_{p}$, is denoted as $\Omega_{i_p}$. Meanwhile, the neighboring cells of 
$\Omega_{i_p}$ are denoted as $\Omega_{i_{pn}}$, where $n=1,\dots, N$ and $N$ is the number of neighboring cells of $\Omega_{i_p}$. 
To achieve the third-order accuracy, a big stencil for non-compact reconstruction for  cell $\Omega_i$  is selected as follows
\begin{align*}
S_i^{WENO}=\{Q_{i},Q_{i_1},\dots, Q_{i_P},Q_{i_{11}},\dots, Q_{i_{1N}},\dots, Q_{i_{P1}},\dots, Q_{i_{PN}}\},
\end{align*}
where $Q_{i_p}$ and $Q_{i_{pn}}$ are the cell averaged conservative variables over cell $\Omega_{i_p}$ and $\Omega_{i_{pn}}$,
respectively. To deal with the discontinuity, the 
sub-stencils $S_{i_m}^{WENO}$ for cell $\Omega_{i}$ are selected, where $m=1,\dots, M$ and $M$ is the number of sub-stencils. 
For the hexahedral cell, the sub-candidate stencils contain the cell averaged conservative variables over cell $\Omega_{i}$ and three neighboring
cells \cite{GKS-WENO-3}. For the tetrahedral cells, to avoid the centroids of $\Omega_{i}$ and three of neighboring cells becoming 
coplanar, the extra cells are added to the sub-candidate stencils \cite{GKS-WENO-3}.  For consistency, 
the cell averaged variable $Q_{i}$ is denoted as $Q_{0}$ and $Q_{m_0}$, and the big stencil and the sub-stencils are rearranged as 
\begin{align*}
S_i^{WENO}&=\{Q_{0},Q_1,…,Q_{K}\},\\
S_{i_m}^{WENO}&=\{Q_{m_0},Q_{m_1},…,Q_{m_K}\},
\end{align*}
where the repeated variables are deleted. 
With the selected stencils, a quadratic polynomial and several linear polynomials can be constructed based on the big stencil $S_i^{WENO}$ 
and the sub-stencils $S_{i_m}^{WENO}$ as follows
\begin{equation}\label{polys}
\begin{split}
P_0(\boldsymbol{x})&=Q_{0}+\sum_{|\boldsymbol d|=1}^2a_{\boldsymbol d}p_{\boldsymbol d}(\boldsymbol{x}),\\
P_m(\boldsymbol{x})&=Q_{m_0}+\sum_{|\boldsymbol d|=1}b_{\boldsymbol d}^mp_{\boldsymbol d}(\boldsymbol{x}),
\end{split}
\end{equation}
where $Q_{i}=Q_{0}=Q_{m_0}$, the multi-index $\boldsymbol d=(d_1, d_2, d_3)$ 
and $|\boldsymbol d|=d_1+d_2+d_3$. The base function $p_{\boldsymbol d}(\boldsymbol{x})$ is defined as
\begin{align}\label{base}
\displaystyle
p_{\boldsymbol d}(\boldsymbol{x})=x^{d_1}y^{d_2}z^{d_3}-\frac{1}{\left|\Omega_{0}\right|}\iiint_{\Omega_{0}}x^{d_1}y^{d_2}z^{d_3}\text{d}V.
\end{align}
To determine these polynomials, the following constrains need to be satisfied
\begin{align*}
\frac{1}{\left|\Omega_{k}\right|}\iiint_{\Omega_{k}}P_0(\boldsymbol{x})\text{d}V&=Q_{k},~Q_{k}\in S_i^{WENO},\\
\frac{1}{\left|\Omega_{m_k}\right|}\iiint_{\Omega_{m_k}}P_m(\boldsymbol{x})\text{d}V&=Q_{m_k},~Q_{m_k}\in S_{i_m}^{WENO}.
\end{align*} 
The over-determined linear systems can be generated, and the least square method is used to obtain the coefficients $a_{\boldsymbol d}$ and $b_{\boldsymbol d}^m$.

\subsection{Compact WENO reconstruction}
For the target cell $\Omega_{i}$, the big stencil of $Q_{i}$ for the compact reconstruction is selected as follows
\begin{align*}
S_i^{HWENO}&=\{Q_{i},Q_{i_1},\dots, Q_{i_P},\nabla Q_{i},\nabla Q_{i_1},\dots,\nabla Q_{i_P}\},
\end{align*}
where $Q_{i_p}$ and $\nabla Q_{i_p}$ are the cell averaged conservative variable and cell averaged gradient of $\Omega_{i_p}$, 
Meanwhile, the sub-stencils $S_{i_m}^{HWENO}$ in the compact HWENO reconstruction for cell $\Omega_{i}$ can be selected more simply. 
The sub-candidate stencils are selected as
\begin{align*}
S_{i_m}^{HWENO}=\{Q_{i},Q_{i_m},\nabla Q_{i},\nabla Q_{i_m}\},
\end{align*}
where $Q_{i_m}$ and $\nabla Q_{i_m}$ are the cell averaged conservative variable and cell averaged gradient over cell $\Omega_{i_m}$, $\Omega_{i_m}$ is one of the neighboring cell of the target cell $\Omega_{i}$, $m=1,\dots,M$, $M$ is the number of sub-stencils which equals to the number of neighboring cells. 
Similar with the WENO reconstruction,  the big stencil and the sub-stencils are rearranged as 
\begin{align*}
S_i^{HWENO}&=\{Q_0,Q_1,…,Q_{P},\nabla Q_0,...,\nabla Q_{P}\},\\
S_{i_m}^{HWENO}&=\{Q_{m_0},Q_{m_1},\nabla Q_{m_0},\nabla Q_{m_1}\}.
\end{align*}
A quadratic polynomial and several linear polynomials in Eq.\eqref{polys} can be constructed based on the big stencil $S_i^{HWENO}$ 
and the sub-stencils $S_{i_m}^{HWENO}$, where the polynomials are defined same as Eq.\eqref{polys} and Eq.\eqref{base}.
To determine these polynomials with a smaller stencil, the additional constrains need to be added for all cells as follows  
\begin{equation*}\label{compact-big-1}
\begin{split}
\frac{1}{\left|\Omega_{k}\right|}\iiint_{\Omega_{k}}P_0(\boldsymbol{x})\text{d}V&=Q_{k},~Q_{k}\in S_i^{HWENO},\\
\frac{h_k}{\left|\Omega_{k}\right|}\iiint_{\Omega_{k}}\frac{\partial}{\partial {\boldsymbol \tau}}P_0(\boldsymbol{x})\text{d}V&=(Q_{\boldsymbol{\tau}})_{k}h_k,~(Q_{\boldsymbol{\tau}})_{k}\in S_i^{HWENO},
\end{split}
\end{equation*}
and
\begin{equation*}\label{compact-big-2}
\begin{split}
\frac{1}{\left|\Omega_{m_k}\right|}\iiint_{\Omega_{m_k}}P_m(\boldsymbol{x})\text{d}V&=Q_{m_k},~Q_{m_k}\in S_{i_m}^{HWENO},\\
\frac{h_{m_k}}{\left|\Omega_{m_k}\right|}\iiint_{\Omega_{m_k}}\frac{\partial}{\partial {\boldsymbol \tau}}P_m(\boldsymbol{z})\text{d}V&
=(Q_{\boldsymbol{\tau}})_{m_k}h_{m_k},~(Q_{\boldsymbol{\tau}})_{m_k}\in S_{i_m}^{HWENO},
\end{split}
\end{equation*}
where $Q_{\boldsymbol{\tau}}$ is the cell averaged derivative, and $\boldsymbol{\tau}$ takes $\boldsymbol{n}_x, \boldsymbol{n}_y, \boldsymbol{n}_z$.
The characteristic lengths $h_k$ and $h_{m_k}$ are used to obtain a small condition number of the matrix to solve $a_{\boldsymbol d}$ and $b_{\boldsymbol d}$. 
The definition of characteristic length is given in Eq.\eqref{length_c}.
The constrained least square method is used to obtain the reconstructed polynomials, where the conservative variable equations are set as strictly satisfied and
others are satisfied in the sense of least square.

In the compact HGKS, the cell averaged derivatives are not updated by the finite volume scheme, which is different from the classical HWENO scheme \cite{HWENO-1,HWENO-2}. 
According to the Gaussian theorem \cite{GKS-HWENO-5}, the cell averaged gradient can be calculated as follows
\begin{equation*}
\begin{split}
|\Omega_k|(\nabla Q)_k(t)&=\iiint_{\Omega_k}\nabla Q(t)\text{d}V=\iint_{\partial\Omega_k}Q(t)\boldsymbol{\tau}\text{d}S=\sum_{i_p\in N(k)} \Big( \sum_{G}\omega_G Q(\boldsymbol{x}_G,t)\boldsymbol{\tau} S_{i_p} \Big),
\end{split}
\end{equation*}
where the conservative variables are obtained by taking moments of gas distribution function 
\begin{align*}
Q(\boldsymbol{x}_G,t)=\int\boldsymbol\psi f(\boldsymbol{x}_G,t,\boldsymbol{u},\xi)\text{d}\Xi.
\end{align*}

\subsection{Non-linear combination}
With the reconstructed polynomial $P_m(\boldsymbol{x}), m=0,...,M$,
the point-value $Q(\boldsymbol{x}_{G})$ and the spatial derivatives
$\partial_{x,y,z} Q(\boldsymbol{x}_{G})$  for reconstructed
variables  at Gaussian quadrature point can be given by the
non-linear combination
\begin{equation*}
\begin{split}
Q(\boldsymbol{x}_{G})=\overline{\omega}_0(\frac{1}{\gamma_0}P_0(\boldsymbol{x}_{G})-&
\sum_{m=1}^{M}\frac{\gamma_m}{\gamma_0}P_m(\boldsymbol{x}_{G}))+\sum_{m=1}^{M}\overline{\omega}_mP_m(\boldsymbol{x}_{G}),\\
\partial_{x,y,z} Q(\boldsymbol{x}_{G})=\overline{\omega}_0(\frac{1}{\gamma_0}\partial_{x,y,z}
P_0(\boldsymbol{x}_{G})-&\sum_{m=1}^{M}\frac{\gamma_m}{\gamma_0}\partial_{x,y,z}
P_m(\boldsymbol{x}_{G}))+\sum_{m=1}^{M}\overline{\omega}_m\partial_{x,y,z}
P_m(\boldsymbol{x}_{G}),
\end{split}
\end{equation*}
where $\gamma_0, \gamma_1,\dots,\gamma_M$ are the linear weights. The
non-linear weights $\omega_m$ and normalized non-linear weights
$\overline{\omega}_m$ are defined as
\begin{align*}
\overline{\omega}_{m}=\frac{\omega_{m}}{\sum_{m=0}^{M} \omega_{m}},~
\omega_{m}=\gamma_{m}\Big[1+\Big(\frac{\tau_Z}{\beta_{m}+\epsilon}\Big)\Big],~
\tau_Z=\sum_{m=1}^{M}\Big(\frac{|\beta_0-\beta_m|}{M}\Big),
\end{align*}
where $\epsilon$ is a small positive number, and $\beta_{m}$ is the smooth indicator. 
More details can be found in \cite{GKS-WENO-3} for 
non-compact reconstruction, and \cite{GKS-HWENO-1, GKS-HWENO-2} for compact reconstruction.

\section{Time-implicit scheme for unsteady flow}

To simulate unsteady flows with complex geometry, the high-order temporal discretizations are used.
Based on the time-dependent flux function of the generalized Riemann problem solver (GRP) \cite{GRP-high-1,GRP-high-2} 
and gas-kinetic scheme (GKS) \cite{GKS-high-1,GKS-high-2},  the two-stage fourth-order time-accurate discretization was developed 
for Lax-Wendroff type flow solvers, particularly applied for the hyperbolic conservation laws.  Considering the semi-discretized finite volume scheme
\begin{align}\label{semi}
\frac{\partial Q_i}{\partial t}=\mathcal{L}(Q_{i}),
\end{align}
where $Q_{i}$ is the cell averaged conservative value of $\Omega_{i}$ and $\mathcal {L}$ is an operator for spatial derivative of flux.

\noindent\textbf{Theorem}:
Introducing an intermediate state at $t^*=t_n+\Delta t/2$, the two-stage temporal discretization for Eq.\eqref{semi} can be written as
\begin{equation}\label{two-stage-1}
\begin{split}
&Q^*_i=Q_i^n+\frac{1}{2}\Delta t\mathcal{L}(Q^n_i)+\frac{1}{8}\Delta
t^2\frac{\partial}{\partial
t}\mathcal{L}(Q^n_i), \\
Q^{n+1}_{i}&=Q^n_i+\Delta t\mathcal {L}(Q^n_i)+\frac{1}{6}\Delta
t^2\big(\frac{\partial}{\partial
t}\mathcal{L}(Q^n_i)+2\frac{\partial}{\partial
t}\mathcal{L}(Q^*_i)\big).
\end{split}
\end{equation}
It can be proved that the two-stage time stepping method Eq.\eqref{two-stage-1} provides a fourth-order time accurate solution for $Q(t)$ at $t=t_n +\Delta t$ \cite{GRP-high-1}.

\subsection{Implicit two-stage third-order temporal discretization}
The two-stage method provides a reliable framework to develop high-order scheme with the implementation of second-order flux function. 
In the explicit two-stage fourth-order GKS, a uniform global time step is used, and the global time step is strongly 
dependent on the smallest cell in the computational mesh. For the unsteady flow with a significant variation of local cell size, 
the global time step can be really small. In order to improve the computational efficiency and keep a high-order accuracy for unsteady 
flows, a two-stage third-order implicit time-accurate discretization is developed.
For explicit two-stage fourth-order GKS, the temporal derivatives $\displaystyle\partial_t\mathcal{L}(Q_i^n)$ and 
$\displaystyle\partial_t\mathcal{L}(Q_i^*)$ in Eq.\eqref{two-stage-1} are determined by solving a linear system \cite{GKS-high-1}. 
For implicit scheme, the temporal derivative $\displaystyle\partial_t\mathcal{L}(Q_i^n)$ in each stage is approximated by 
the backward Euler scheme. In the first stage of Eq.\eqref{two-stage-1}, $\displaystyle\partial_t\mathcal{L}(Q_i^n)$ is approximated by
\begin{align*}
\displaystyle\frac{\partial}{\partial t}\mathcal{L}(Q_i^n)&=\frac{\mathcal{L}(Q_i^*)-\mathcal{L}(Q_i^n)}{\Delta t/2}.
\end{align*}
In the second stage of Eq.\eqref{two-stage-1}, $\displaystyle\partial_t\mathcal{L}(Q_i^n)$ and 
$\displaystyle\partial_t\mathcal{L}(Q_i^*)$ are given by
\begin{align*}
\displaystyle\frac{\partial}{\partial t}\mathcal{L}(Q_i^n)&=\frac{\mathcal{L}(Q_i^{n+1})-\mathcal{L}(Q_i^n)}{\Delta t},~
\displaystyle\frac{\partial}{\partial t}\mathcal{L}(Q_i^*)=\frac{\mathcal{L}(Q_i^*)-\mathcal{L}(Q_i^n)}{\Delta t/2}.
\end{align*}

\noindent\textbf{Theorem}:
Thus, the implicit two-stage temporal discretization for Eq.\eqref{semi} can be written as
follows
\begin{equation}\label{two-stage-2}
\begin{split}
&Q_i^*=Q_i^n+ \frac{\Delta t}{4}\left(\mathcal{L}(Q_i^n)+\mathcal{L}(Q_i^*)\right),\\
Q_i^{n+1}&=Q_i^n+\frac{\Delta t}{6}(\mathcal{L}(Q_i^{n})+4\mathcal{L}(Q_i^{*})+\mathcal{L}(Q_i^{n+1})).
\end{split}
\end{equation}
It can be proved that the two-stage time stepping method Eq.\eqref{two-stage-2} provides a third-order time accurate solution for $Q(t)$ at $t=t_n +\Delta t$.

\noindent\textbf{Proof}:
Integrating Eq.\eqref{semi} on the time interval $[t^n,t^{n+1}]$, we have
\begin{align*}
Q_i^{n+1}-Q_i^n&=\int_{t_n}^{t_n+\Delta t}\mathcal{L}(Q_i(t))\text{d}t.
\end{align*}
To prove that Eq.\eqref{two-stage-2} provides a third-order accurate approximation, the following Taylor expansion need to be satisfied  
\begin{align}\label{exp}
\int_{t_n}^{t_n+\Delta t}\mathcal{L}(Q(t))\text{d}t=\Delta t\mathcal{L}(Q^n)+\frac{\Delta t^2}{2}\mathcal{L}_t(Q^n)
+\frac{\Delta t^3}{6}\mathcal{L}_{tt}(Q^n)+\mathcal{O}(\Delta t^4).
\end{align}
According to the Cauchy-Kovalevskaya method, the temporal derivatives are given by
\begin{align*}
&\mathcal{L}_{t}=\mathcal{L}_Q\mathcal{L},\\
&\mathcal{L}_{tt}=\mathcal{L}_Q^2\mathcal{L}+\mathcal{L}_{QQ}\mathcal{L}^2,\\
&\mathcal{L}_{ttt}=\mathcal{L}_Q^3\mathcal{L}+4\mathcal{L}_{QQ}\mathcal{L}_{Q}\mathcal{L}^2+\mathcal{L}_{QQQ}\mathcal{L}^3.
\end{align*}
For the operator $\mathcal{L}$, we have the following expansion up
to the corresponding order
\begin{align*} 
\mathcal{L}(Q^*)=\mathcal{L}(Q^n)&+\mathcal{L}_{Q}(Q^*-Q^n)+\frac{\mathcal{L}_{QQ}}{2}(Q^*-Q^n)^2\\
&+\frac{\mathcal{L}_{QQQ}}{6}(Q^{*}-Q^n)^3+\mathcal{O}(Q^*-Q^n)^4,
\end{align*}
and
\begin{align*}
\mathcal{L}(Q^{n+1})=\mathcal{L}(Q^n)&+\mathcal{L}_Q(Q^{n+1}-Q^n)+\frac{\mathcal{L}_{QQ}}{2}(Q^{n+1}-Q^n)^2\\
                     &+\frac{\mathcal{L}_{QQQ}}{6}(Q^{n+1}-Q^n)^3+\mathcal{O}(Q^{n+1}-Q^n)^4.
\end{align*}
Substituting above expansion into Eq.\eqref{two-stage-2}, we have
\begin{equation*}\label{proof2}
\begin{split}
Q^*-Q^n&=\frac{\Delta t}{2}\mathcal{L}(Q^n)+\frac{\Delta t}{4}\mathcal{L}_Q(Q^*-Q^n)+\frac{\Delta t}{8}\mathcal{L}_{QQ}(Q^*-Q^n)^2+\Delta t \mathcal{O}(Q^*-Q^n)^3\\
&=\frac{\Delta t}{2}\mathcal{L}(Q^n)+\frac{\Delta t}{4}\mathcal{L}_Q(\frac{\Delta t}{2}\mathcal{L}(Q^n)+\frac{\Delta t^2}{8}\mathcal{L}_t(Q^n))+\frac{\Delta t}{8}\mathcal{L}_{QQ}(\frac{\Delta t^2}{4}\mathcal{L}^2(Q^n))+O(\Delta t^4)\\
&=\frac{\Delta t}{2}\mathcal{L}(Q^n)+\frac{\Delta t^2}{8}\mathcal{L}_t(Q^n)+\frac{\Delta t^3}{32}\mathcal{L}_{tt}(Q^n)+\mathcal{O}(\Delta t^4)
\end{split}
\end{equation*}
and
\begin{equation*}\label{proof3}
\begin{split}
Q^{n+1}-Q^n&=\Delta t\mathcal{L}(Q^n)+\frac{\Delta t}{6}\big(\mathcal{L}_Q(Q^{n+1}-Q^n)+4\mathcal{L}_Q(Q^*-Q^n)\big)\\
&\quad\quad+\frac{\Delta t}{12}\big(\mathcal{L}_{QQ}(Q^{n+1}-Q^n)^2+4\mathcal{L}_{QQ}(Q^*-Q^n)^2\big)+\mathcal{O}(\Delta t^4)\\
&=\Delta t\mathcal{L}(Q^n)+\frac{\Delta t^2}{2}\mathcal{L}_t(Q^n)+\frac{\Delta t^3}{6}\mathcal{L}_{tt}(Q^n)+\mathcal{O}(\Delta t^4),
\end{split}
\end{equation*}
which is identical to the right hand side of Eq.\eqref{exp}.
Thus, it has be proved that Eq.\eqref{two-stage-2} provides the third-order temporal accuracy for Eq.\eqref{semi} at $t=t_n+\Delta t$.

\subsection{Implicit two-stage HGKS for unsteady flows}
For the unsteady Euler and Navier-Stokes equations, the dual time-stepping method was adopted for the implicit scheme \cite{GMRES-4}.  
Eq.\eqref{two-stage-2} can be rewritten into the following form
\begin{equation}\label{two-stage-3}
\begin{split}
\frac{Q_{i}^*-Q_{i}^n}{\Delta t/2}&=\frac{1}{2}\left(\mathcal{L}(Q_{i}^n)+\mathcal{L}(Q_{i}^*)\right),\\
\frac{Q_{i}^{n+1}-Q_{i}^n}{\Delta t}&=\frac{1}{6}(\mathcal{L}(Q_{i}^{n})+4\mathcal{L}(Q_{i}^{*})+\mathcal{L}(Q_{i}^{n+1})).
\end{split}
\end{equation} 
Each stage in Eq.\eqref{two-stage-3} can be considered as a nonlinear system, which needs to be solved. 
In order to solve the systems with the pseudo-time iteration, the artificial temporal derivative term is added into each stage
\begin{equation}\label{artificial}
\begin{split}
&\frac{Q_{i}^{m+1}-Q_{i}^m}{\Delta t_a}=-\frac{Q_{i}^{m+1}-Q_{i}^n}{\Delta t/2}+\frac{1}{2}\left(\mathcal{L}(Q_{i}^n)+\mathcal{L}(Q_{i}^{m+1})\right),\\
&\frac{Q_{i}^{m+1}-Q_{i}^m}{\Delta t_a}=-\frac{Q_{i}^{m+1}-Q_{i}^n}{\Delta t}+\frac{1}{6}(\mathcal{L}(Q_{i}^{n})+4\mathcal{L}(Q_{i}^{*})+\mathcal{L}(Q_{i}^{m+1})),
\end{split}
\end{equation}
where $\Delta t_a$ is  the artificial time step, $t^{m+1}=t^{m}+\Delta t_a$, $t^{m=0}=t^n$. In each stage, $Q_{i}^{m}$ converges 
to $Q_{i}^*$ and $Q_{i}^{n+1}$.  $\mathcal{L}(Q_{i}^{m+1})$ in each stage can be linearized as
\begin{equation*}
\mathcal{L}(Q_{i}^{m+1})=\mathcal{L}(Q_{i}^{m})+\displaystyle(\frac{\rm{d} \mathcal{L}}{\rm{d} Q})_i^m\Delta Q_{i}^m.
\end{equation*}
Substituting the linearization above into artificial iterations, Eq.\eqref{artificial} can be written as
\begin{equation}\label{two-stage-4}
\begin{split}
&\displaystyle\big[(\frac{1}{\Delta t_a}+\frac{1}{\Delta t/2})\boldsymbol{I}-\frac{1}{2}(\frac{\rm{d} \mathcal{L}}{\rm{d} Q})_i^m\big]\Delta Q_{i}^m=\frac{Q_{i}^{n}-Q_{i}^m}{\Delta t/2}+\frac{1}{2}\left(\mathcal{L}(Q_{i}^n)+\mathcal{L}(Q_{i}^{m})\right),\\
&\displaystyle\big[(\frac{1}{\Delta t_a}+\frac{1}{\Delta t})\boldsymbol{I}-\frac{1}{6}(\frac{\rm{d} \mathcal{L}}{\rm{d} Q})_i^m\big]\Delta Q_{i}^m=\frac{Q_{i}^{n}-Q_{i}^m}{\Delta t}+\frac{1}{6}(\mathcal{L}(Q_{i}^{n})+4\mathcal{L}(Q_{i}^{*})+\mathcal{L}(Q_{i}^{m})),
\end{split}
\end{equation}
where $\boldsymbol{I}$ is the identity matrix. Eq.\eqref{two-stage-4} can be solved by the classical implicit methods. 
In this paper, both the lower-upper symmetric Gauss-Seidel (LUSGS) method  \cite{LUSGS-2}  
and the nonlinear generalized minimal residual (GMRES) method  \cite{GMRES-1,GMRES-2}  are adopted.

\begin{algorithm}[!h]
\SetAlgoLined
Initial condition for $t^n$\;
\While{\rm{TIME} $ \leq $ \rm{TSTOP}}{
    Calculation of time step $\Delta t$\;
    WENO/HWENO reconstruction to calculate $\mathcal{L}(Q^n)$\;
    \textcolor{blue}{STAGE 1:}\\
    Initial condition for $t^m$\;
    \For{$i=1, k_a $}{
    Calculation of time step $\Delta t_a$\;
    WENO/HWENO reconstruction to calculate $\mathcal{L}(Q^m)$\;
    LUSGS/GMRES for solving artificial steady stage 1\;    
    Artificial iteration: $Q^{m+1}=Q^{m}+\Delta Q^m$, $Q_{\boldsymbol{\tau}}^{m+1}$\;}
    Update 1: $Q^{*}=Q^{m+1}$, $Q_{\boldsymbol{\tau}}^*=Q_{\boldsymbol{\tau}}^{m+1}$\;
    WENO/HWENO reconstruction to calculate $\mathcal{L}(Q^*)$\;
    \textcolor{blue}{STAGE 2:}\\
    Initial condition for $t^m$\;
    \For{$i=1, k_a $}{
    Calculation of time step $\Delta t_a$\;
    WENO/HWENO reconstruction to calculate $\mathcal{L}(Q^m)$\;
    LUSGS/GMRES for solving artificial steady stage 2\;    
    Artificial iteration: $Q^{m+1}=Q^{m}+\Delta Q^m$, $Q_{\boldsymbol{\tau}}^{m+1}$\;}
    Update 2: $Q^{n+1}=Q^{m+1}$, $Q_{\boldsymbol{\tau}}^{n+1}=Q_{\boldsymbol{\tau}}^{m+1}$\;
}
\caption{\label{S2O3-algorithm} Program for two-stage third-order implicit HGKS scheme}
\end{algorithm}

For the LUSGS method, the differentiation in the left side of Eq.\eqref{two-stage-4} can be approximated by the Euler equations-based flux splitting method as follows
\begin{equation*}
\displaystyle(\frac{\rm{d} \mathcal{L}}{\rm{d} Q})_i\Delta Q_{i}^m=\mathcal{L}(Q_{i}^{m+1})-\mathcal{L}(Q_{i}^{m})=-\frac{1}{|\Omega_{i}|}\sum_{i_p\in N(i)}\Delta F_{i,i_p}^{m}S_{i_p},
\end{equation*}
where
\begin{equation*}
\Delta F_{i,i_p}^{m}=\frac{1}{2}(\Delta T_{i_p}^m+\Delta T_i^m-\lambda_{i_p}(\Delta Q_{i_p}^m-\Delta Q_i^m)).
\end{equation*}
In the equation above, $\Delta T^m=T^{m+1}-T^m$ and $\Delta Q^m=Q^{m+1}-Q^m$,
$T_{i}$ and $T_{i_p}$ are the Euler flux corresponding to cell $\Omega_i$ and  $\Omega_{i_p}$, 
and the spectral radius of the Euler flux Jacobian $\lambda_{i,i_p}$ satisfies the following inequality
\begin{equation*}
\lambda_{i,i_p}\geq|\boldsymbol{u}_{i,i_p}\cdot \boldsymbol{n}_{i,i_p}|+a_{i,i_p},
\end{equation*}
where $\boldsymbol{u}_{i,i_p}$ and $a_{i,i_p}$ are the velocity and the speed of sound at
the interface $S_{i_p}$ and $\boldsymbol{n}_{i,i_p}$ is the unit normal direction of $S_{i_p}$.
With the approximation, the Gauss-Seidel iteration process can be applied to solve Eq.\eqref{two-stage-4} by a forward sweep step and a backward sweep step.

For the GMRES method, according to the total differential formulation, the differentiation in the left side of Eq.\eqref{two-stage-4} can be written as
\begin{equation*}
(\frac{\rm{d} \mathcal{L}}{\rm{d} Q})_i\Delta Q_{i}^m=-\frac{1}{|\Omega_{i}|}\sum_{i_p\in N(i)}(\frac{\partial F}{\partial Q_i})^m\Delta Q_i^mS_{i_p}
-\frac{1}{|\Omega_{i}|}\sum_{i_p\in N(i)}(\frac{\partial F}{\partial Q_{i_p}})^m\Delta Q_{i_p}^mS_{i_p}
\end{equation*}
The Jacobian matrix about the cell $\Omega_i$ and  $\Omega_{i_p}$ can be derived by
\begin{equation*}
\begin{cases}
\displaystyle\frac{\partial F}{\partial Q_i}=\frac{1}{2}(J(Q_i)+\lambda_{i,i_p}I),\\
\displaystyle\frac{\partial F}{\partial Q_{i_p}}=\frac{1}{2}(J(Q_{i_p})-\lambda_{i,i_p}I),
\end{cases}
\end{equation*}
where $J(Q)$ represents the Jacobian matrix of the inviscid flux and $\lambda_{i,i_p}$ is same as the spectral radius in LUSGS method.
With the approximation, Eq.\eqref{two-stage-4} can be solved by GMRES method as a large sparse nonlinear system.
In the computation, the Jacobi iteration is used as the preconditioner.

For the compact schemes, the cell averaged directional derivatives also need to be updated at each time step.
Due to the implicit methods are applied in the pseudo steady problems, the update of cell averaged directional derivatives can be driven by time evolution directly.
In the artificial iterations,
\begin{align*}
(Q_{\boldsymbol{\tau}})_i^{n+1}&=
\frac{1}{|\Omega_i|}\sum_{i_p\in N(i)}\Big(\sum_{G}\omega_G\Big(\int\psi f(\boldsymbol{x}_G,t^n,\boldsymbol{u},\xi)\boldsymbol{\tau}\text{d}\Xi\Big)S_{i_p}\Big).
\end{align*}
where $\boldsymbol{\tau}=\boldsymbol{n}_x,\boldsymbol{n}_y,\boldsymbol{n}_z$.
Algorithm.\ref{S2O3-algorithm} shows the whole process of two-stage third-order implicit HGKS scheme, where $k_a$ is the artificial iteration times.

\noindent\textbf{Remrak}:
Usually, the serial sweep part in LUSGS method is not easy to be implemented for parallel computation.
Meanwhile, for GMRES method, the parallel computation can be implemented, except the minimization process.
Although GMRES method is more complicated, it is expected that it would be more efficient than LUSGS method.

\section{Numerical tests}
In this section, the numerical tests for both inviscid and viscous flows  will be presented to validate the current schemes. 
For the sake of convenience, the implicit two-stage third-order HGKS combined with LUSGS method and GMRES method are denoted as 
S2O3-L and S2O3-G, respectively. Meanwhile, the explicit two-stage fourth-order HGKS is denoted as S2O4-E.
In the computation, the dimension of Krylov subspace is  denoted as $\dim K$, which is set as 3 without special statement.
The times of artificial iterations is denoted as $k_a$, which will be mentioned for each case.

For evolution of flow fields, the time step $\Delta t$ is given by the CFL condition
\begin{equation*}
\Delta t=\text{CFL}\cdot\mathop{min}\limits_{i}\frac{h_i}{(|\boldsymbol{u}_i|+a_i)},
\end{equation*}
where $|\boldsymbol{u}_i|$ is the magnitude of velocity, $a_i$ is the speed of sound and the characteristic length of $\Omega_i$ is defined as 
\begin{equation}\label{length_c}
h_i=\displaystyle\frac{|\Omega_i|}{\max\limits_{i_p} |S_{i_p}|},
\end{equation}
where $|S_{i_p}|$ is the area of $S_{i_p}$. For iteration of artificial steady problems,
the artificial time step $\Delta t_a$ is denoted by artificial CFL number $\text{CFL}_\text{a}$. For all test cases, $\text{CFL}_\text{a}$ is set as 2000.
In addition, the time step $\Delta t_s$ is introduced for evaluating the time-averaged flux. For explicit schemes, $\Delta t_s=\Delta t$. 
For implicit schemes, $\Delta t_s$ is given by $\text{CFL}\leq 1$. The time step $\Delta t_s$ is also used for calculating collision time.
For the inviscid flows, the collision time $\tau$ takes
\begin{align*}
\tau=\epsilon \Delta t_s+C\displaystyle|\frac{p_l-p_r}{p_l+p_r}|\Delta t_s,
\end{align*}
where $p_l$ and $p_r$ denote the pressure on the left and right sides of the cell interface, $\epsilon=0.05$ and $C=2.5$. For the viscous flows, we have
\begin{align*}
\tau=\frac{\mu}{p}+C \displaystyle|\frac{p_l-p_r}{p_l+p_r}|\Delta t_s,
\end{align*}
where $\mu$ is the dynamic viscous coefficient and $p$ is the pressure at the cell interface. In smooth
flow regions, the collision time reduces to
\begin{equation*}
\tau=\frac{\mu}{p}.
\end{equation*}

\begin{figure}[!h]
\centering
\includegraphics[width=0.4\textwidth]{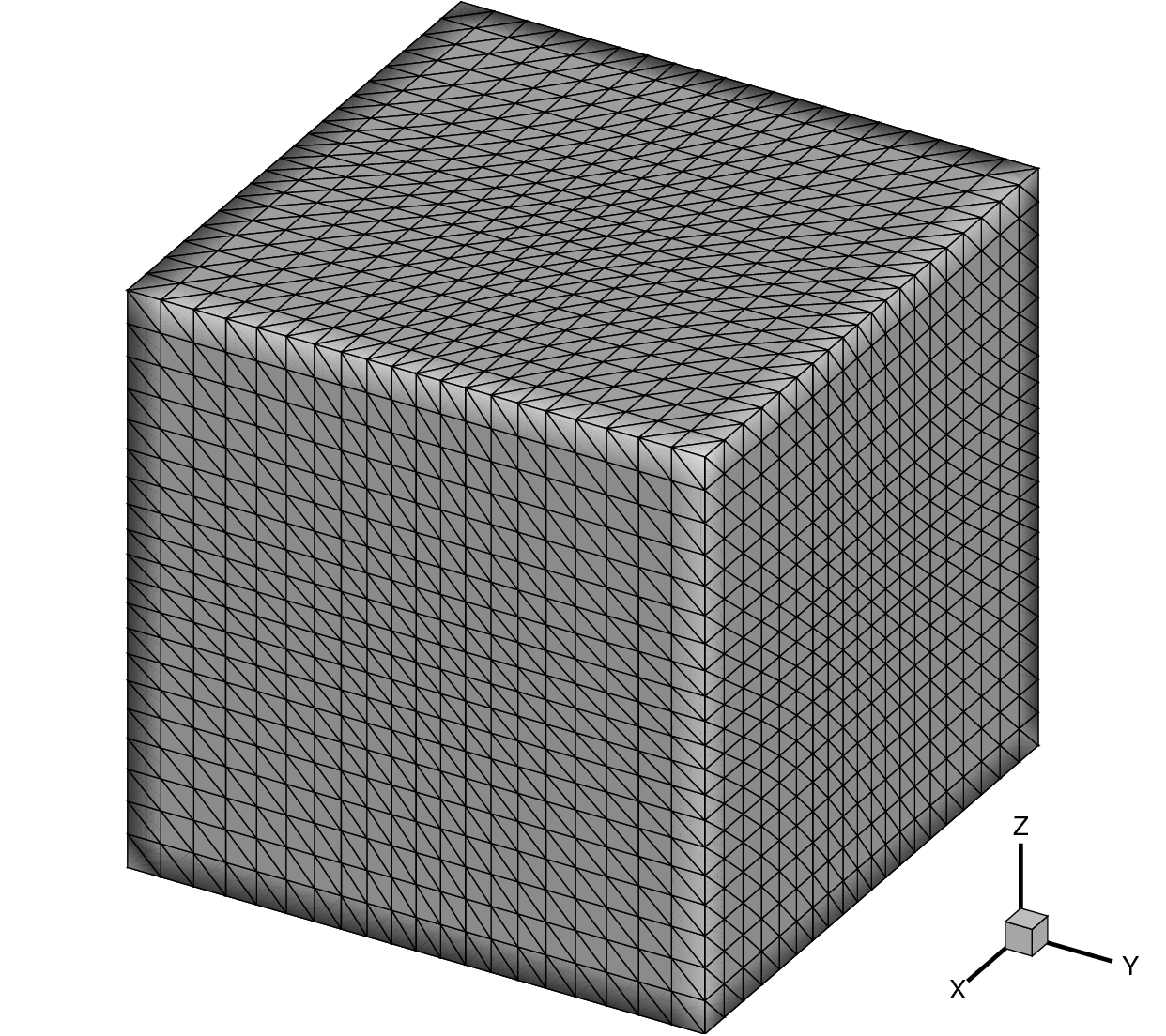}
\caption{\label{accuracy-mesh} Accuracy test: the computational mesh distribution with $20^3\times6$ tetrahedral cells.}
\end{figure}

\begin{table}[!h]
\begin{center}
\def\temptablewidth{0.93\textwidth}{\rule{\temptablewidth}{1.0pt}}
\begin{tabular*}{\temptablewidth}{@{\extracolsep{\fill}}c|c|c|c|c|c} 
MESH            & noncompact S2O4-E & Order &  MESH            & compact S2O4-E & Order  \\
\hline
$5^3\times 6$   & $1.6310\times 10^{-1}$ &     &
$5^3\times 6$   & $6.5181\times 10^{-2}$ &         \\
$10^3\times 6$  & $2.8220\times 10^{-2}$ & 2.5310 &
$10^3\times 6$  & $8.6376\times 10^{-3}$ & 2.9157  \\
$20^3\times 6$  & $3.7253\times 10^{-3}$ & 2.9213 &
$20^3\times 6$  & $1.0528\times 10^{-3}$ & 3.0364  \\
$40^3\times 6$  & $4.7009\times 10^{-4}$ & 2.9863 &
$40^3\times 6$  & $1.3018\times 10^{-4}$ & 3.0157  
\end{tabular*}
{\rule{\temptablewidth}{1.0pt}}
\end{center}
\caption{\label{tab-3d-1} Accuracy tests: the errors and orders of accuracy for the explicit S2O4-E method with tetrahedron meshes.}
\begin{center}
\def\temptablewidth{0.93\textwidth}{\rule{\temptablewidth}{1.0pt}}
\begin{tabular*}{\temptablewidth}{@{\extracolsep{\fill}}c|c|c|c|c|c} 
MESH            & noncompact S2O3-L & Order &  MESH            & compact S2O3-L & Order  \\
\hline
$5^3\times 6$   & $1.9326\times 10^{-1}$ &     &
$5^3\times 6$   & $7.9342\times 10^{-2}$ &         \\
$10^3\times 6$  & $3.7102\times 10^{-2}$ & 2.3810 &
$10^3\times 6$  & $1.2896\times 10^{-2}$ & 2.6212  \\
$20^3\times 6$  & $4.9486\times 10^{-3}$ & 2.9064 &
$20^3\times 6$  & $1.2713\times 10^{-3}$ & 3.3426  \\
$40^3\times 6$  & $6.2522\times 10^{-4}$ & 2.9846 &
$40^3\times 6$  & $1.5226\times 10^{-4}$ & 3.0617  
\end{tabular*}
{\rule{\temptablewidth}{1.0pt}}
\end{center}
\caption{\label{tab-3d-2} Accuracy tests: the errors and orders of accuracy for the implicit S2O3-L method with  tetrahedron meshes.}
\begin{center}
\def\temptablewidth{0.93\textwidth}{\rule{\temptablewidth}{1.0pt}}
\begin{tabular*}{\temptablewidth}{@{\extracolsep{\fill}}c|c|c|c|c|c} 
MESH            & noncompact S2O3-G & Order &  MESH            & compact S2O3-G & Order  \\
\hline
$5^3\times 6$   & $1.9326\times 10^{-1}$ &     &
$5^3\times 6$   & $7.9336\times 10^{-2}$ &         \\
$10^3\times 6$  & $3.7103\times 10^{-2}$ & 2.3810 &
$10^3\times 6$  & $1.2880\times 10^{-2}$ & 2.6228 \\
$20^3\times 6$  & $4.9488\times 10^{-3}$ & 2.9064 &
$20^3\times 6$  & $1.2660\times 10^{-3}$ & 3.3468 \\
$40^3\times 6$  & $6.2532\times 10^{-4}$ & 2.9844 & 
$40^3\times 6$  & $1.5047\times 10^{-4}$ & 3.0727 
\end{tabular*}
{\rule{\temptablewidth}{1.0pt}}
\end{center}
\caption{\label{tab-3d-3} Accuracy tests: the errors and orders of accuracy for the implicit S2O3-G method with  tetrahedron meshes.}
\end{table}

\subsection{Accuracy test}
In this case, the three-dimensional advection of density perturbation is used 
to test the order of accuracy of the non-compact and compact implicit HGKSs. 
The computational domain is $[0,2]\times[0,2]\times[0,2]$. A series of tetrahedron 
meshes with $6 \times N^3$ cells are used, where every cubic is divided into six tetrahedron cells. 
The mesh distribution is shown in Figure.\ref{accuracy-mesh}. The initial condition is given as follows
\begin{align*}
\rho_0&(x, y, z)=1+0.2\sin(\pi(x+y+z)),~p_0(x,y,z)=1,\\
&U_0(x,y,z)=1,~V_0(x,y,z)=1,~W_0(x,y,z)=1.
\end{align*}
The periodic boundary condition is applied on all boundaries, and the
exact solution is
\begin{align*}
\rho(x,y&,z,t)=1+0.2\sin(\pi(x+y+z-3t)),~p(x,y,z,t)=1,\\
&U(x,y,z,t)=1,~V(x,y,z,t)=1,~W(x,y,z,t)=1.
\end{align*}
In this case,  $\tau=0$ and the gas-distribution function reduces to
\begin{align*}
f(\boldsymbol{x}_{G},t,\boldsymbol{u},\xi)=g(1+At).
\end{align*}
To avoid affecting the computation of errors, the residuals of the artificial steady problems must converge to a sufficiently small size. 
The times of artificial iteration takes $k_a=5$ for both S2O3-L and S2O3-G methods and $\dim K=10$ for S2O3-G method. 
The non-compact WENO and compact HWENO schemes with linear reconstructions are tested, respectively. 
The $L^2$ errors and orders of accuracy obtained by both explicit and implicit methods $t=2$ are presented in 
Table.\ref{tab-3d-1}, Table.\ref{tab-3d-2} and Table.\ref{tab-3d-3}. 
It can be observed that the expected third-order of accuracy is achieved by the implicit S2O3-L and S2O3-G methods. 
Based on the same temporal discretization method, 
the compact HGKSs have smaller errors than the non-compact HGKSs due to the small stencils for reconstruction.

\begin{figure}[!h]
\centering
\includegraphics[width=0.45\textwidth]{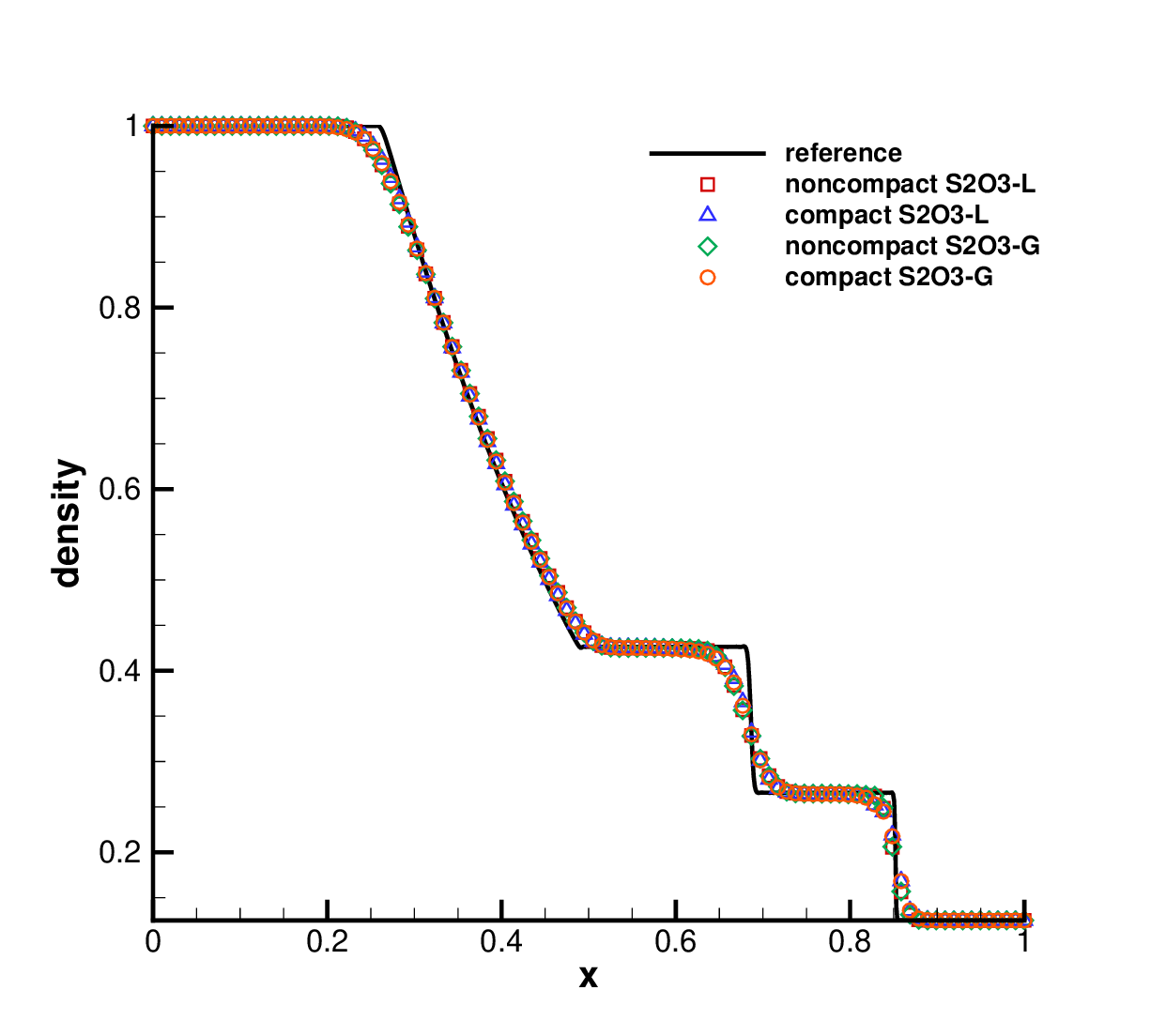} 
\includegraphics[width=0.45\textwidth]{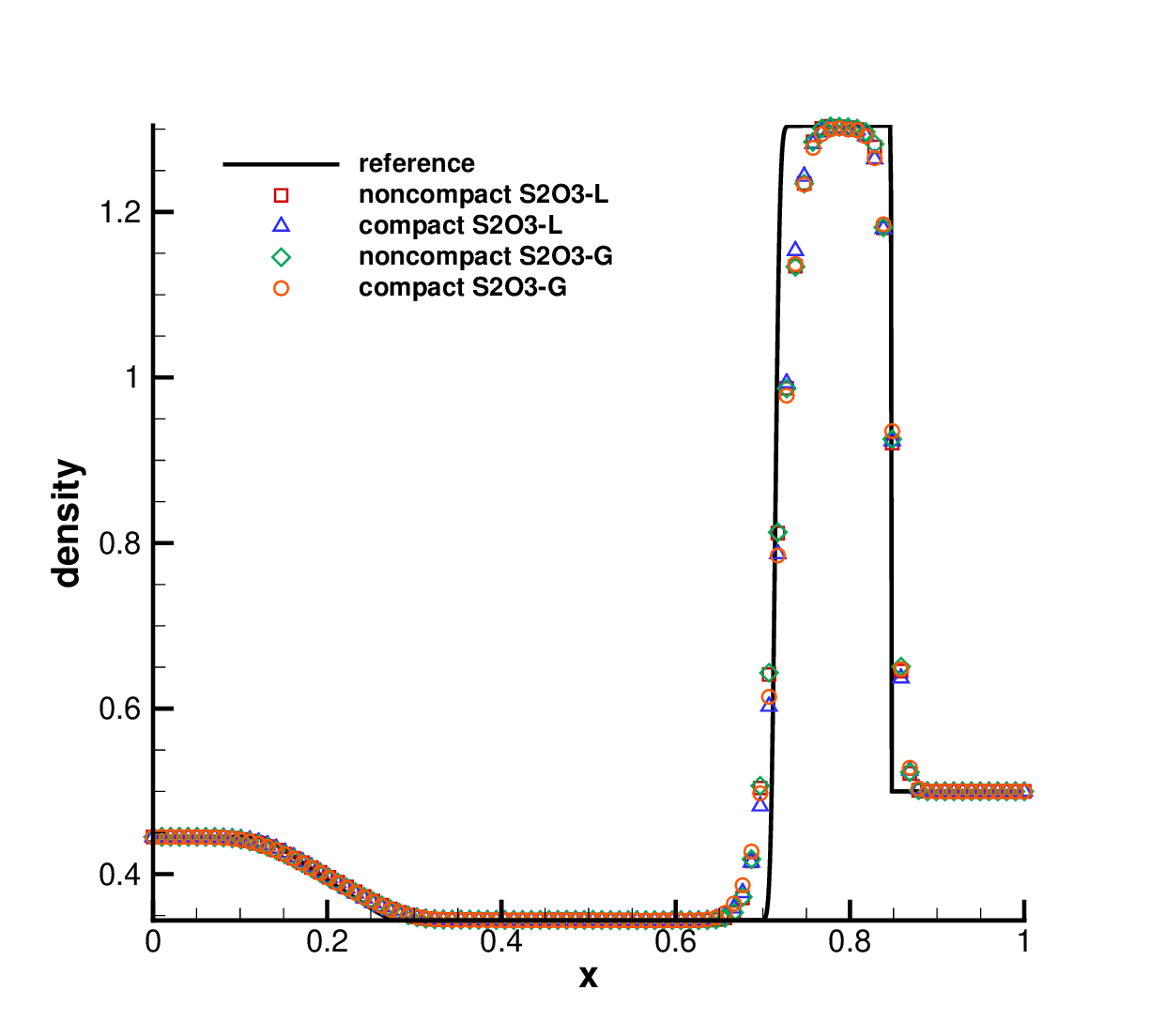}\\
\includegraphics[width=0.45\textwidth]{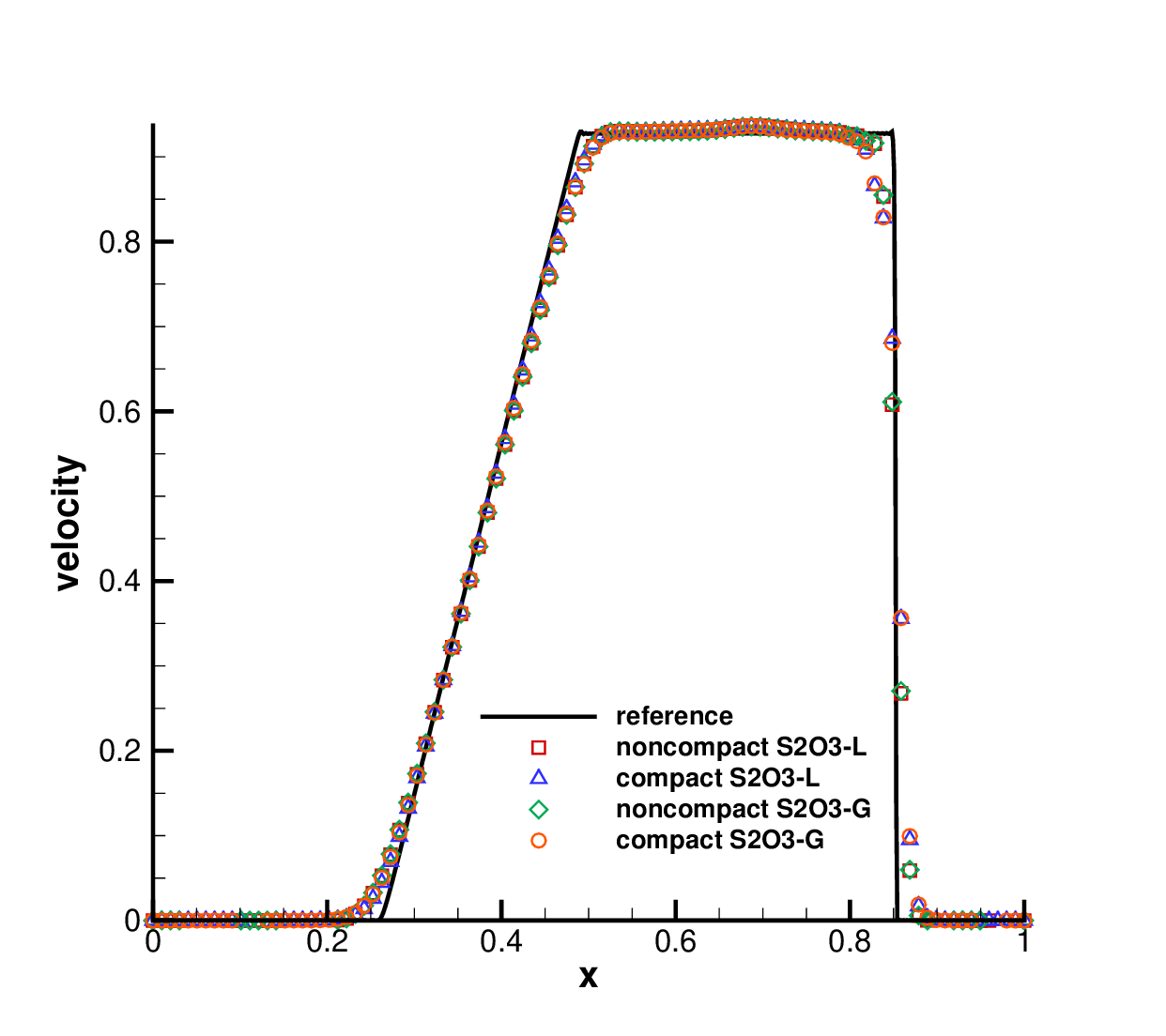} 
\includegraphics[width=0.45\textwidth]{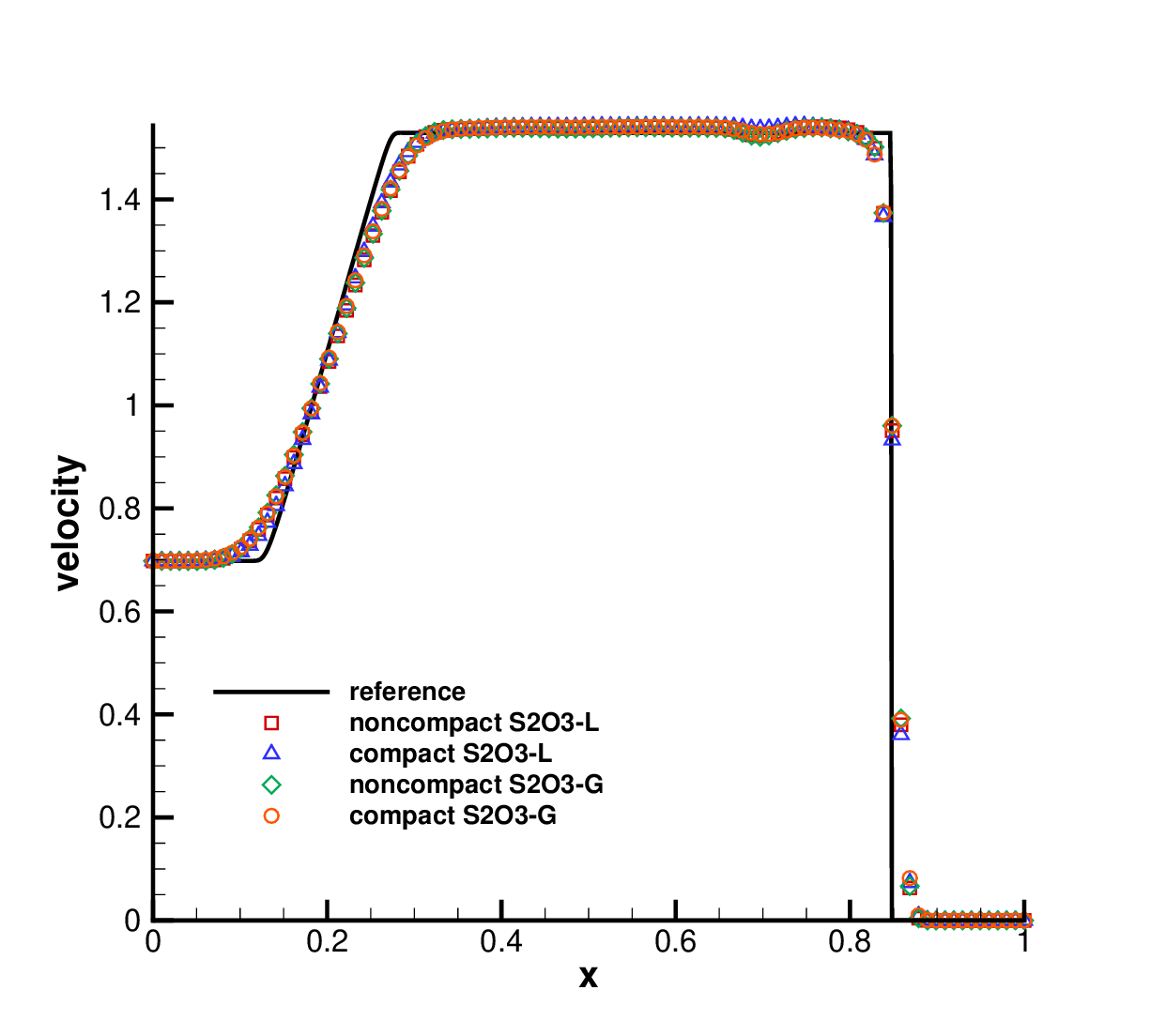}\\
\includegraphics[width=0.45\textwidth]{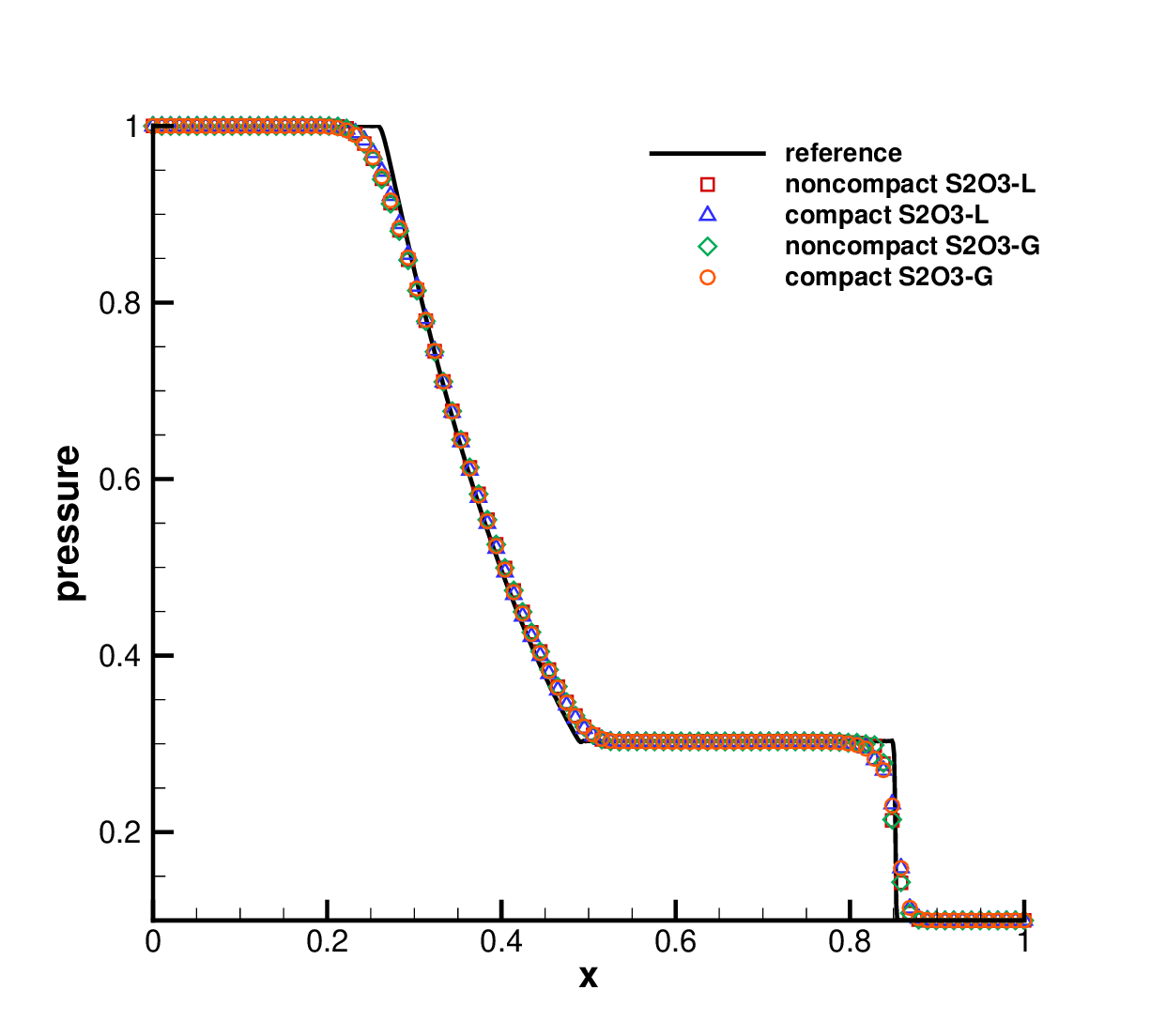} 
\includegraphics[width=0.45\textwidth]{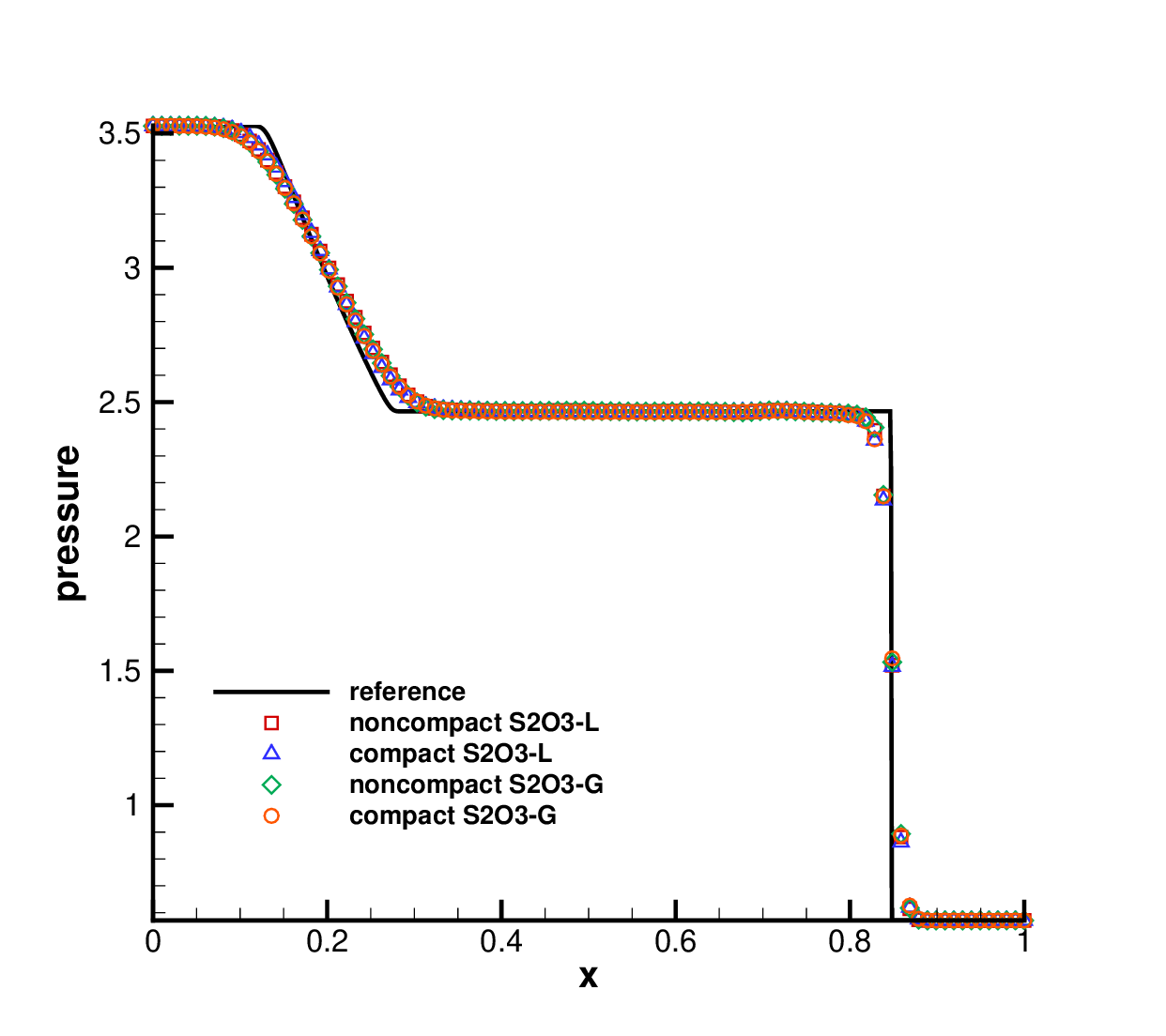}
\caption{\label{Riemann-1d} One-dimensional Riemann problem: the density, velocity and pressure 
distributions for Sod problem (left) at $t=0.2$ and Lax problem (right) at $t=0.14$ at the center horizontal line.}
\end{figure}

\subsection{One-dimensional  Riemann problem}
In this case, the one-dimensional Riemann problems are provided to validate the performance of implicit scheme with shocks.
The non-compact and compact S2O3-L and S2O3-G methods are tested with the three-dimensional hexahedral mesh.
The computational domain takes $[0,1]\times[0,0.1]\times[0,0.1]$. The uniform hexahedral mesh with $100\times 5\times 5$ 
cells is used. The first one is the Sod problem, and the initial condition is given as follows
\begin{equation*}
(\rho,U,V,W,p) = \begin{cases}
(1,0,0,0,1),  &0\leq x<0.5,\\
(0.125, 0,0,0, 0.1), &0.5\leq x\leq1.
\end{cases}
\end{equation*}
The second one is the Lax problem, and the initial condition is given as follows
\begin{equation*}
(\rho,U,V,W,p) = \begin{cases}
(0.445, 0.698,0,0, 3.528),   &0\leq x<0.5,\\
(0.5, 0, 0,0,0.571), &0.5\leq x\leq1.
\end{cases}
\end{equation*}
Non-reflection boundary condition is adopted at all boundaries of the computational domain.
The times of artificial iterations $k_a$ is set as 4 for the Sod problem and 3 for the Lax problem.
For the Sod problem, the CFL number is set as 2.5 for all the four implicit methods. 
Under this condition, it only takes 17 steps to complete the computation, where the 
smallest time step in the whole process satisfies $\Delta t > 1.1\times 10^{-2}$. 
For the Lax problem, the CFL number is set as 3.0 and it takes 22 steps to complete the 
computation, where the smallest time step satisfies $\Delta t > 6.3\times 10^{-3}$.
The density, velocity and pressure distributions for the Sod problem at $t=0.2$ and the 
Lax problem at $t=0.14$ are presented in Figure.\ref{Riemann-1d} with $x\in[0,1]$ and $y=z=0.05$. 
The numerical results agree well with the exact solutions. In the precious implicit schemes \cite{GKS-implicit-unsteady-2,GKS-implicit-unsteady-3,GMRES-4}, 
the spurious oscillations usually occur with the Crank-Nicolson temporal discretization. For the current schemes, the discontinuities 
are well resolved without spurious oscillations, even with such big time steps.

\begin{figure}[!h]
\centering
\includegraphics[width=0.4\textwidth]{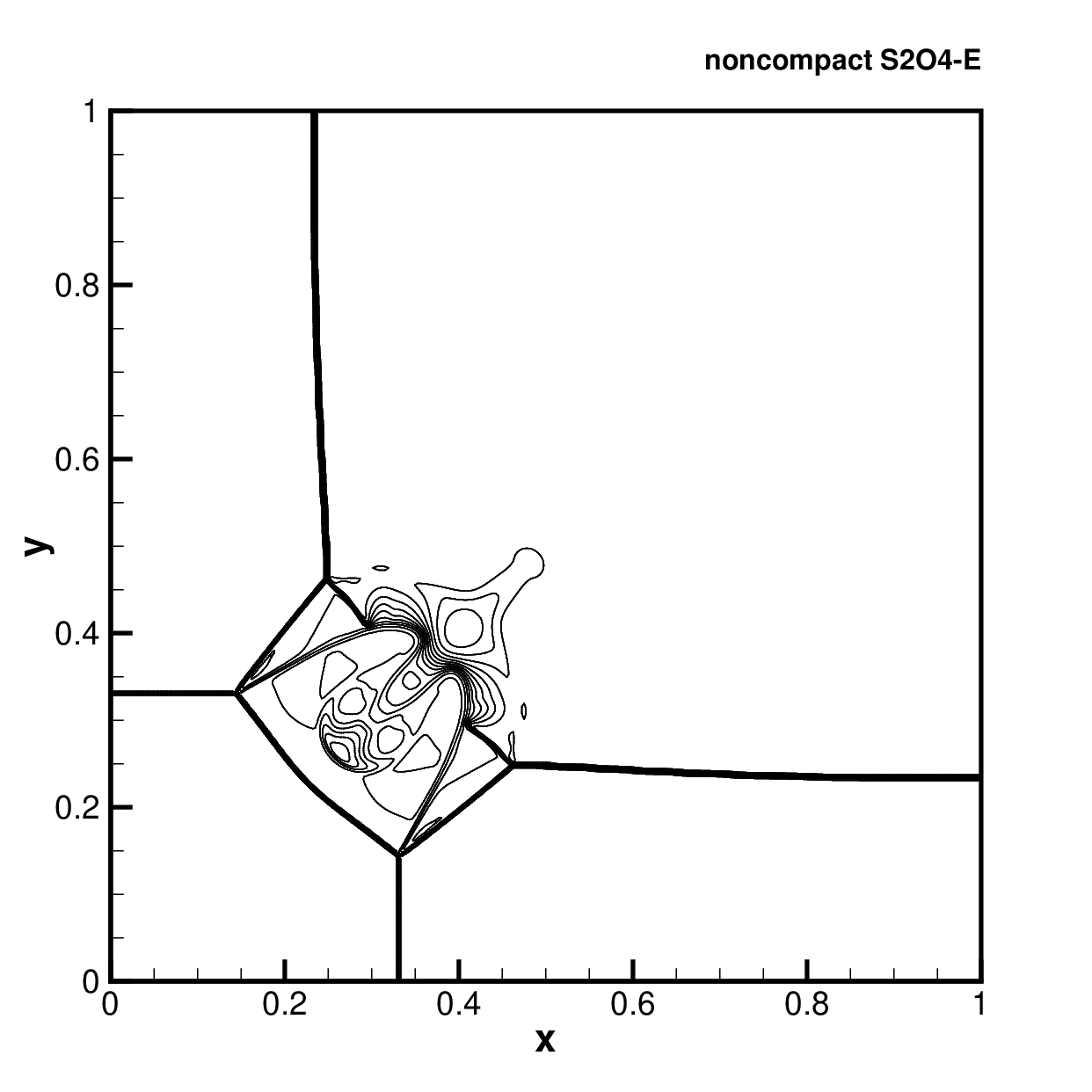}
\includegraphics[width=0.4\textwidth]{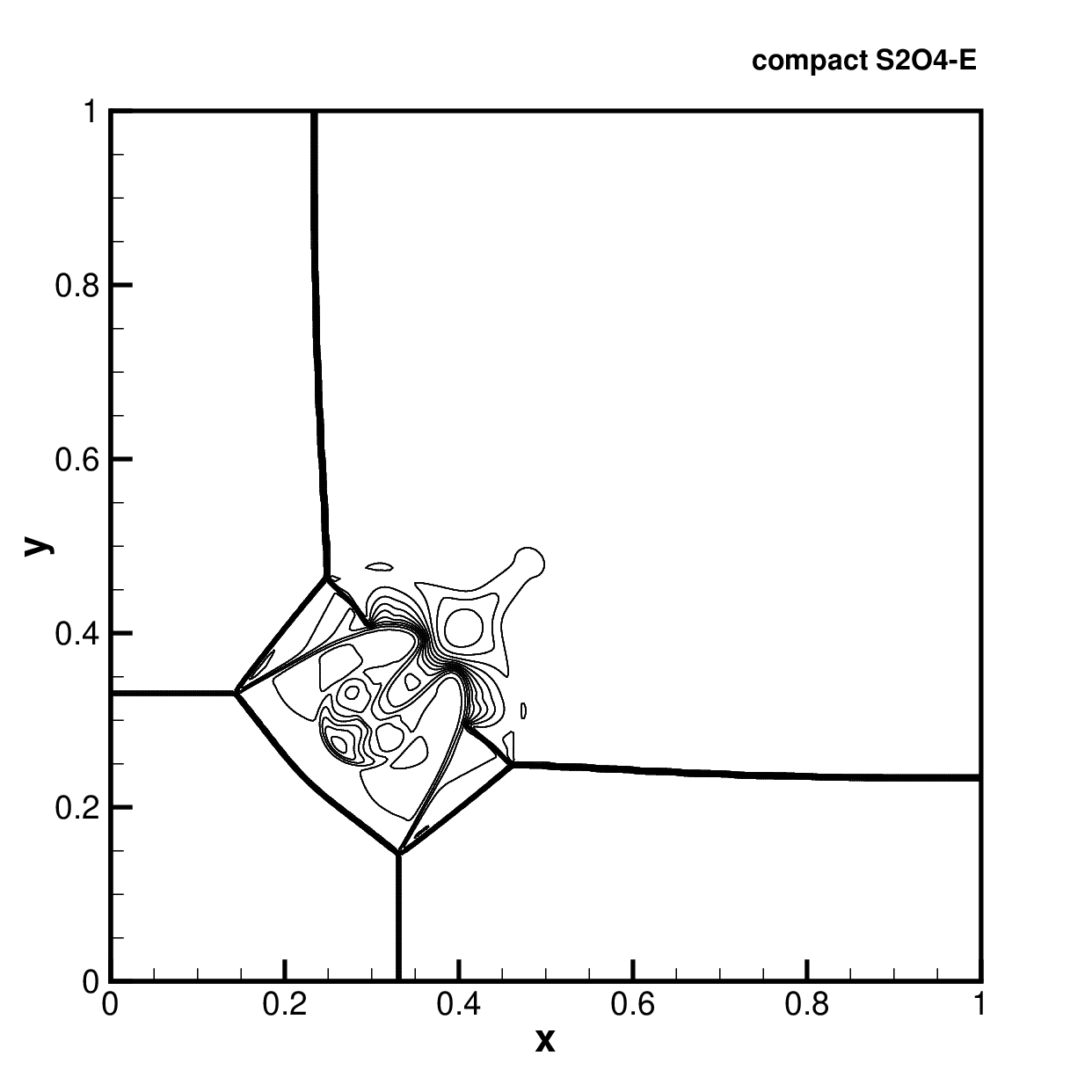}\\
\includegraphics[width=0.4\textwidth]{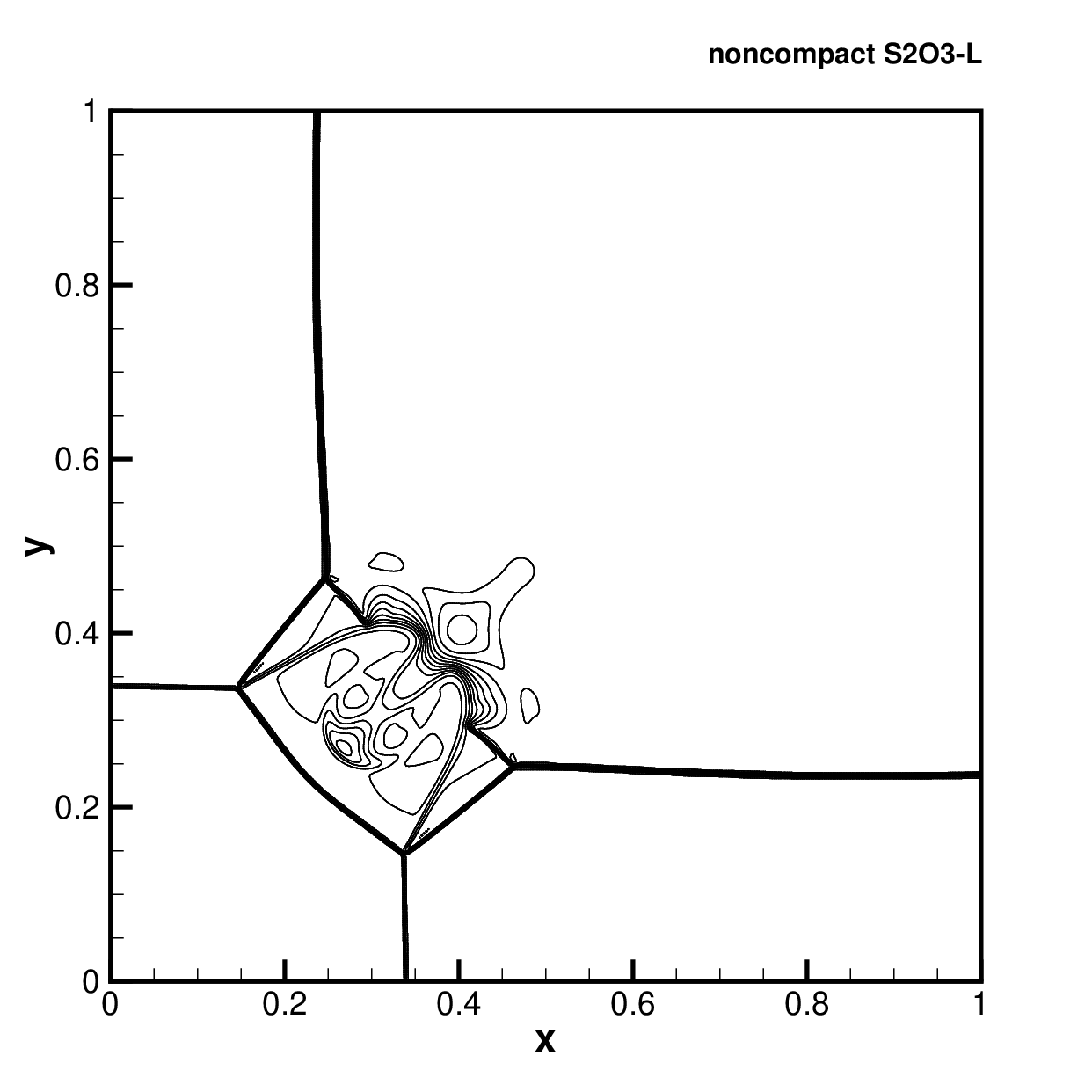}
\includegraphics[width=0.4\textwidth]{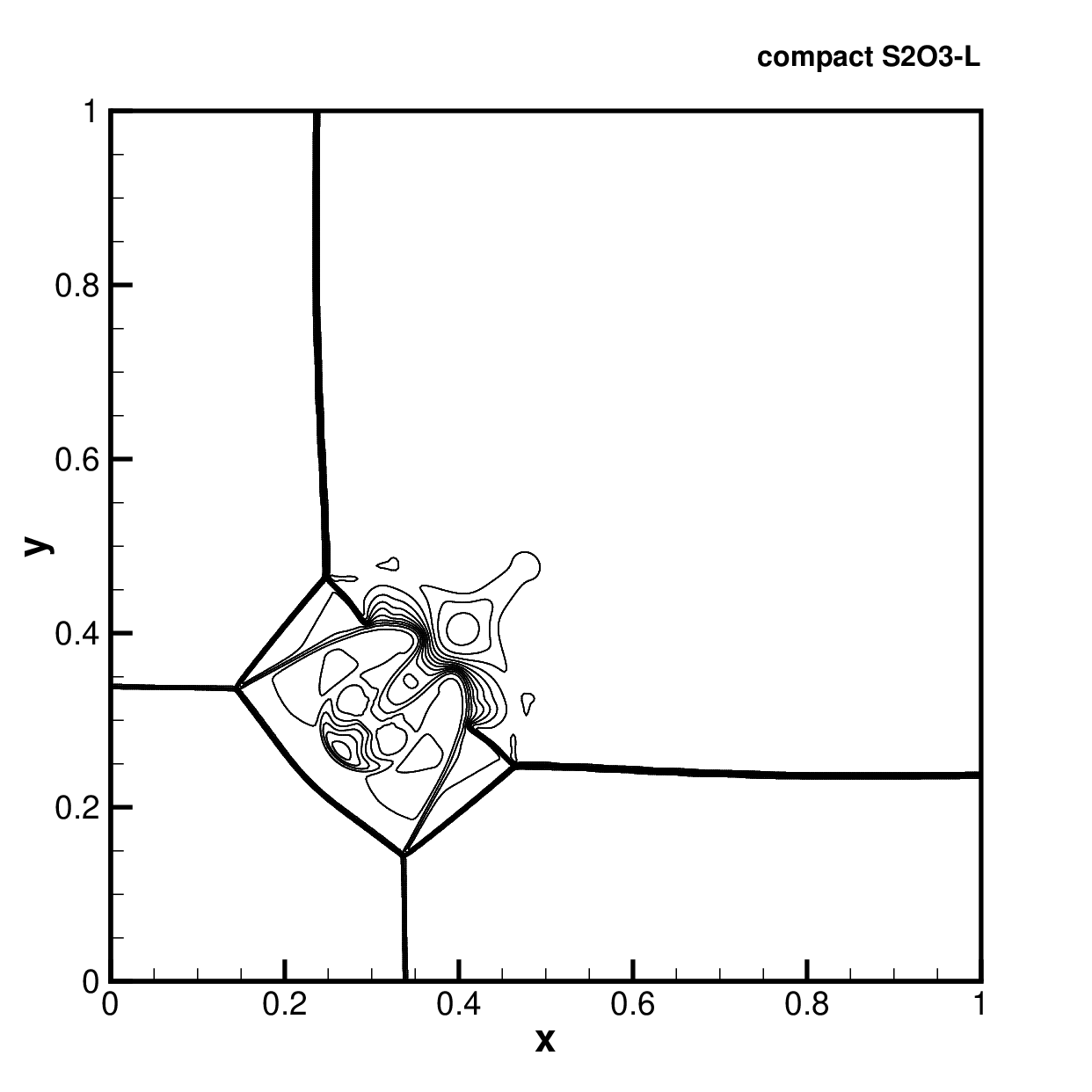}\\
\includegraphics[width=0.4\textwidth]{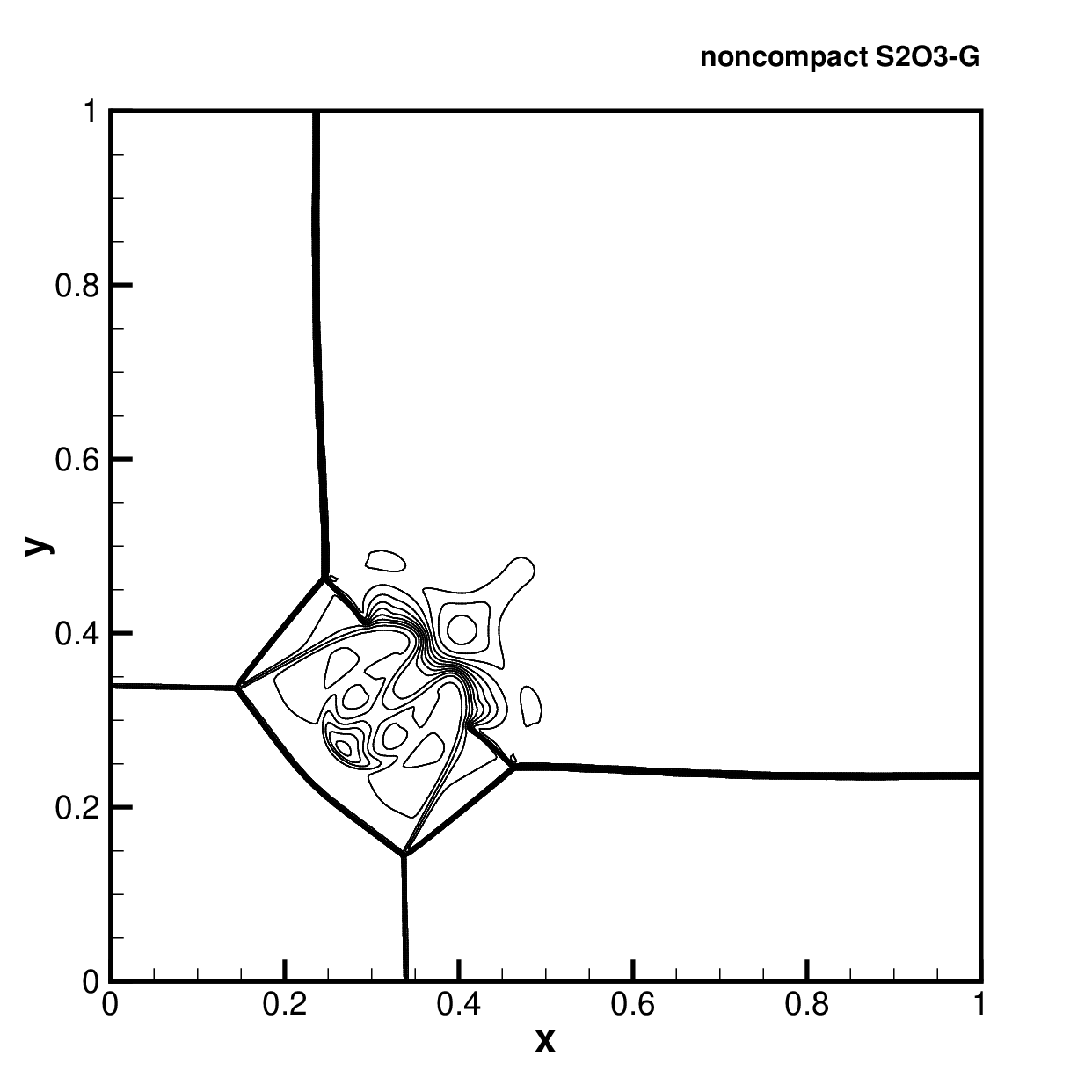}
\includegraphics[width=0.4\textwidth]{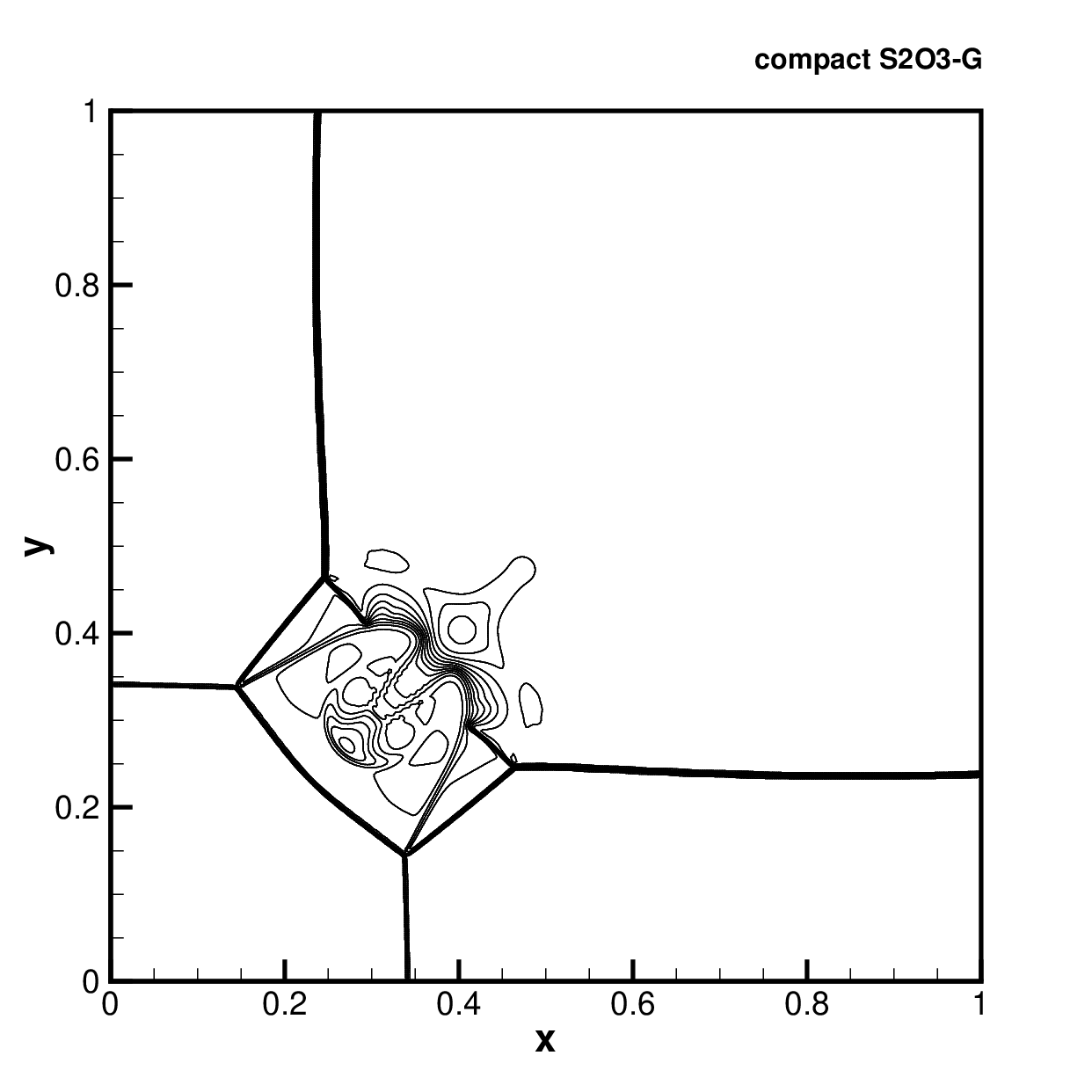}
\caption{\label{Riemann-2d} Two-dimensional  Riemann problem: the density distribution of S2O4-E (top), 
S2O3-L (middle) and S2O3-G (bottom) methods $t=0.4$.}
\end{figure}

\subsection{Two-dimensional Riemann problem}
In this case, the two-dimensional Riemann problem for Euler equations is presented \cite{Case-Lax}, 
and the interaction of four shock waves is provided to validate the robustness of implicit schemes.
The computational domain is $[0,1]\times[0,1]\times[0,0.03]$. 
The uniform hexahedral mesh with $400\times 400\times 3$ cells is used, and 
the non-reflection condition is used for all boundaries. 
The initial condition is given as follows
\begin{equation*}
(\rho,U,V,W,p) = \begin{cases}
(1.5,0,0,0,1.5),  x>0.5, y>0.5,\\
(0.5323, 1.206, 0, 0,0.3), x<0.5, y>0.5,\\
(0.138, 1.206, 1.206, 0, 0.029), x<0.5, y<0.5,\\
(0.5323, 0, 1.206,0, 0.3), x>0.5, y<0.5.
\end{cases}
\end{equation*}
The initial shock wave $S^-_{23}$ bifurcates at the trip point into a reflected shock wave, a Mach stem, and a slip line.
The reflected shock wave interacts with the shock wave $S^-_{12}$ to produce a new shock. 
In the computation, the CFL number is set as 0.6 for explicit methods and 3.0 for implicit methods, and 
the times of artificial iterations takes $k_a=3$. For the non-compact and compact S2O3-L and S2O3-G methods, 
the smallest time steps in the whole computation satisfy $\Delta t > 2.5\times 10^{-3}$. 
Meanwhile, for the non-compact and compact S2O4-E methods, the smallest time steps $\Delta t \approx 5.0\times 10^{-4}$.
The density distributions of S2O4-E, S2O3-L and S2O3-G with non-compact and compact reconstruction
are given in Figure.\ref{Riemann-2d} at $t = 0.4$. The current implicit schemes well resolve the the structures of shock waves, 
and the results indicate the robustness of implicit methods for strong discontinuities.

\begin{figure}[!h]
\centering
\includegraphics[width=0.675\textwidth]{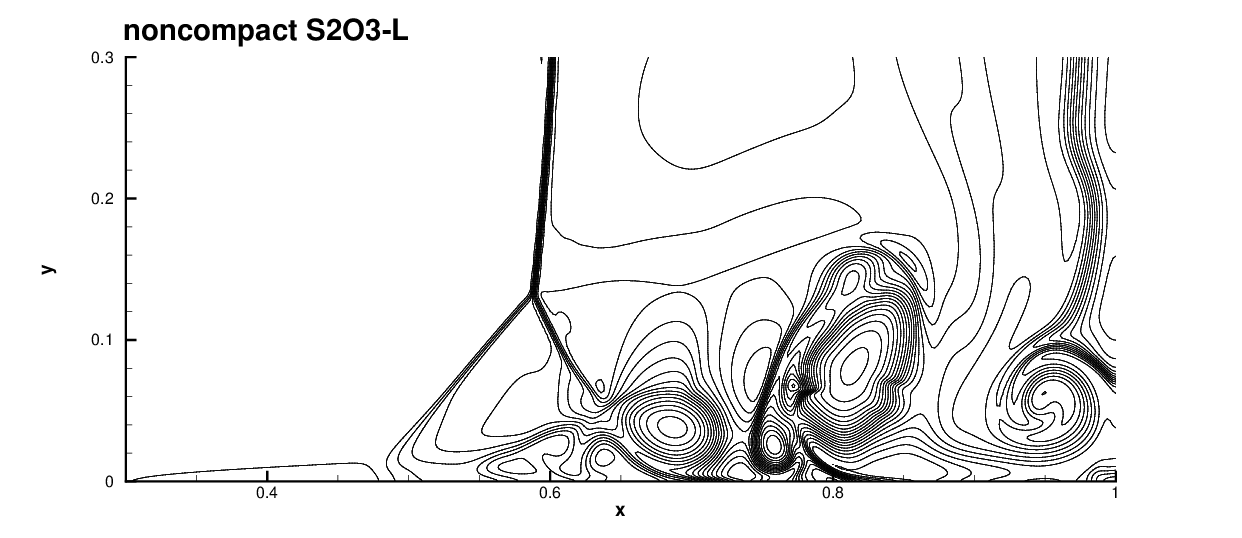}
\includegraphics[width=0.675\textwidth]{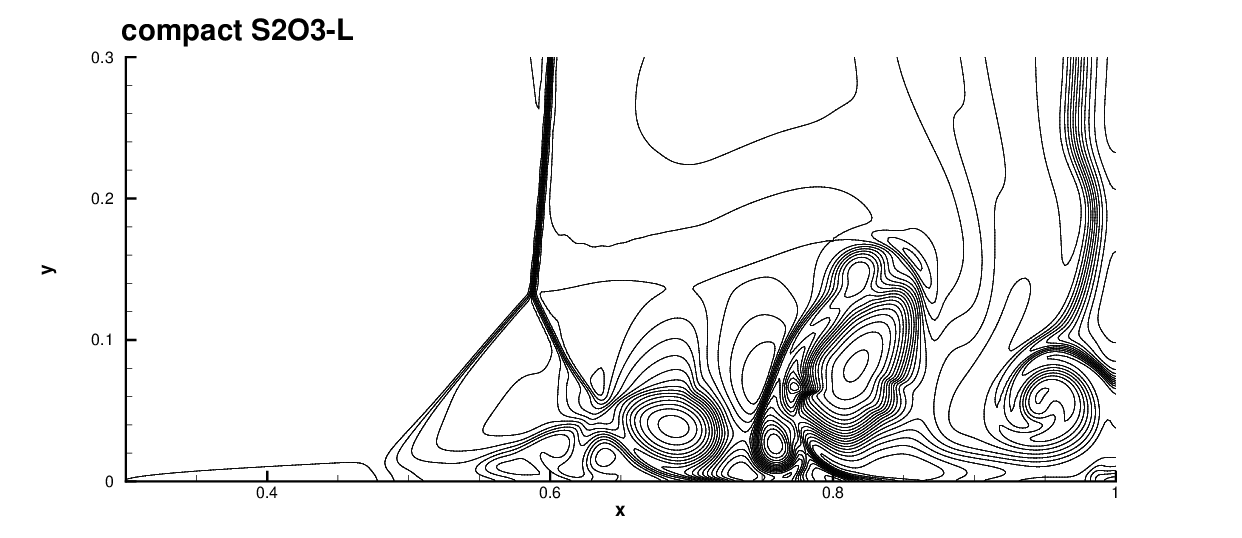}
\caption{\label{Shock-tube-1a} Viscous shock tube: the density distributions with non-compact and compact S2O3-L at $t=1.0$.}
\centering
\includegraphics[width=0.675\textwidth]{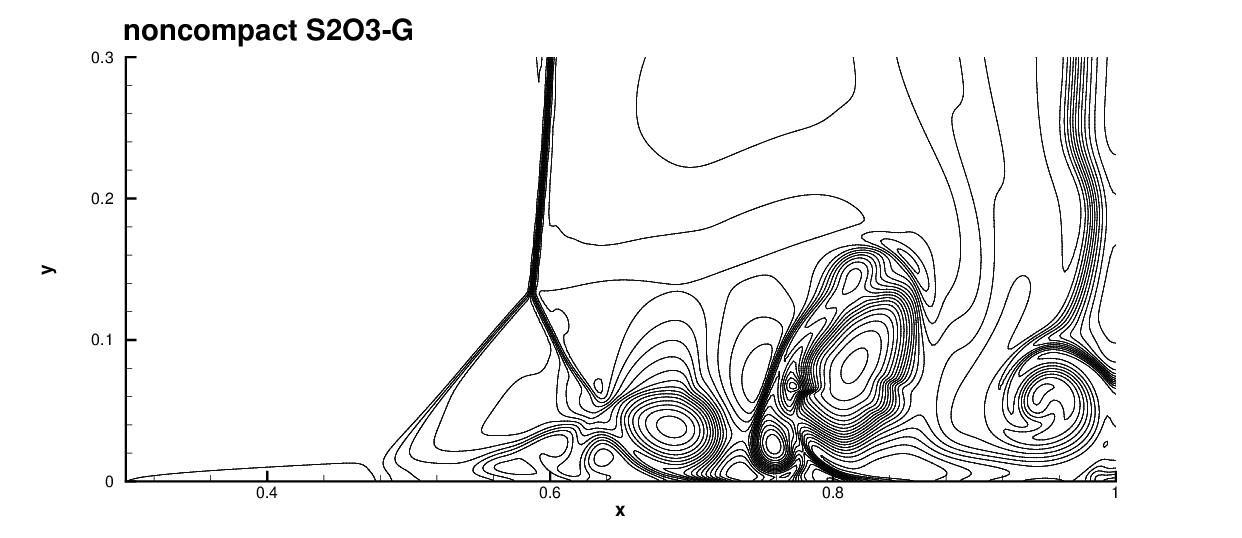}
\includegraphics[width=0.675\textwidth]{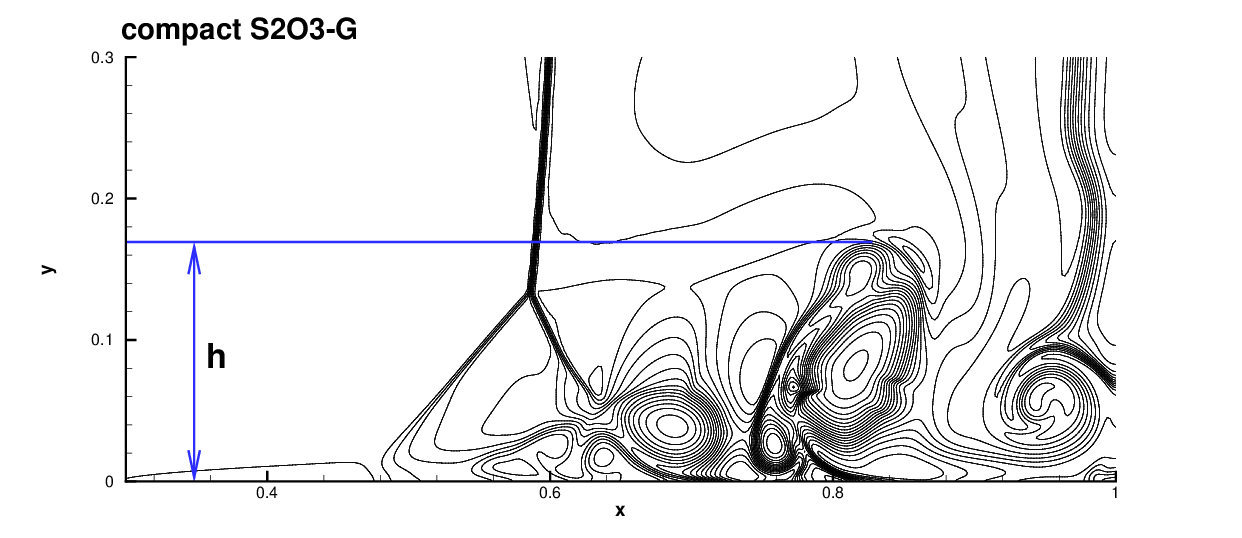}
\caption{\label{Shock-tube-1b} Viscous shock tube: the density distributions with non-compact and compact S2O3-G at $t=1.0$.}
\end{figure}

\begin{figure}[!h]
\centering
\includegraphics[width=0.7\textwidth]{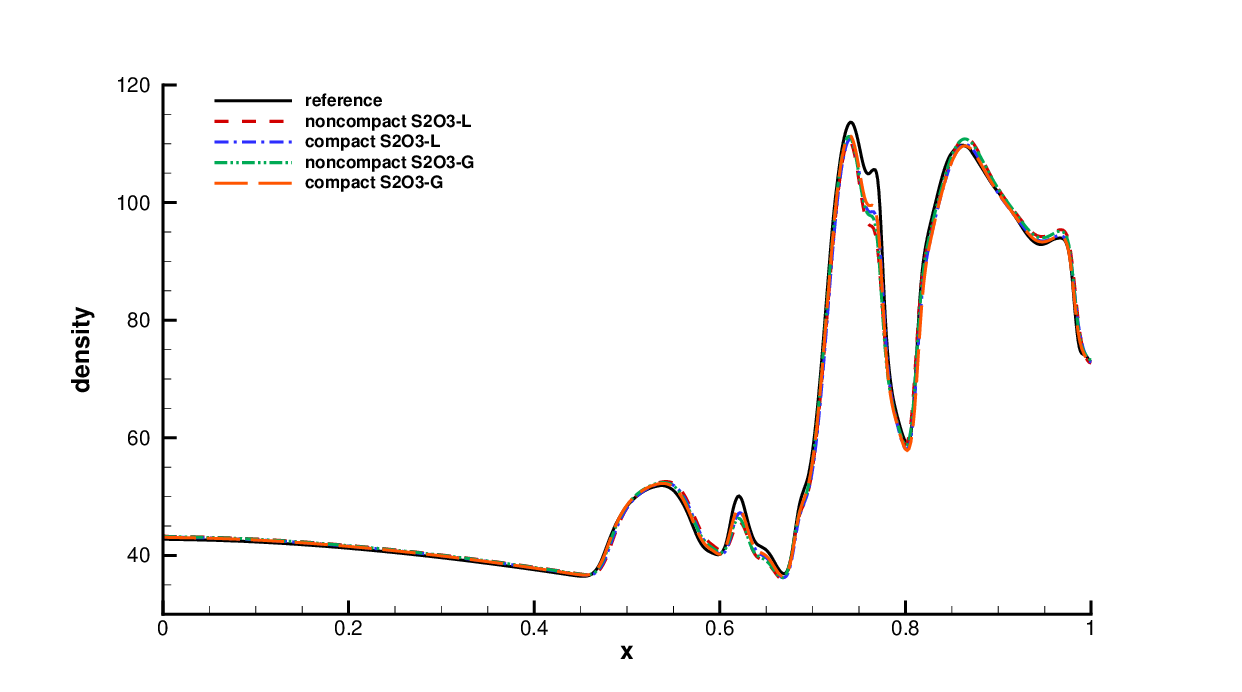}
\caption{\label{Shock-tube-2} Viscous shock tube: the density distribution along the lower wall at $t=1.0$.}
\end{figure}

\subsection{Viscous shock tube}
The viscous shock tube problem  \cite{Case-Daru} is  introduced to test the performances of the current schemes for viscous flows.
The ideal gas is at rest in a two-dimensional unit box $[0, 1]\times [0, 1]$, and the membrane 
located at $x = 0.5$ separates two different states of the gas.  
The initial condition is given by
\begin{equation*}
(\rho,U,V,W,p) = \begin{cases}
(120,0,0,0,120/\gamma),  0<x<0.5,\\
(1.2,0,0,0,1.2/\gamma), 0.5<x<1.0,
\end{cases}
\end{equation*}
where $\gamma=1.4$ and $Re=200$. The membrane is removed at time zero and wave interaction occurs. 
A shock wave moves to the right and reflects at the right end wall. 
After the reflection, it interacts with the contact discontinuity. The contact discontinuity and shock wave 
interact with the horizontal wall and create a thin boundary layer during their propagation.

\begin{table}[!h]
\begin{center}
\begin{tabular}{cc}
\toprule
Scheme                                                        & height \\
\midrule
Fourth-order HGKS  \cite{GKS-high-1}        & 0.171   \\
M-AUSMPW+           \cite{Case-Kim}           & 0.168  \\
AUSMPW+               \cite{Case-Kim}           &  0.163  \\
non-compact S2O3-L  & 0.166   \\  
compact S2O3-L         & 0.169  \\  
non-compact S2O3-G  &  0.167 \\  
compact S2O3-G         & 0.170\\
\bottomrule
\end{tabular}
\caption{\label{Shock-tube-3} Comparison of the heights of primary vortex among current schemes and  
reference data for the reflecting shock-boundary layer interaction.}
\end{center}
\end{table}

Symmetric boundary condition is used on the top boundary $y = 0.5$ and the $z$-direction boundaries. 
Non-slip boundary condition for velocity, and adiabatic condition for temperature are imposed at solid wall boundaries.
In the actual computation, this case is tested in the domain $[0, 1]\times[0, 0.5]\times [0,0.006]$, and the hexahedral mesh which contains $500\times 250\times 3$ cells is used.
The mesh near the wall $y=0$ is refined by the following formulation 
\begin{equation*}
\begin{cases}
x=x,\\
y=y-\sin(2\pi y)/12.5,\\
z=z.
\end{cases}
\end{equation*}
The CFL number is set as 3.5 and the times of artificial iterations $k_a$ is set as 3 for 
non-compact and compact S2O3-L and S2O3-G methods. The density distributions at 
$t = 1.0$ are given in Figure.\ref{Shock-tube-1a} and Figure.\ref{Shock-tube-1b} for implicit S2O3-L and S2O3-G methods. 
The numerical results computed by different implicit methods match well with each other.
The density profiles along the lower wall with $x\in [0,1]$ and $z=0.003$ are given in Figure.\ref{Shock-tube-2}.
The implicit results agree well with the reference data,  where the reference data is calculated by fourth-order HGKS on two-dimensional 
uniform structured mesh with $750\times 750$ cells \cite{GKS-high-1}.
Besides, the height of primary vortex is compared to verify the resolution of current schemes.
As shown in Table.\ref{Shock-tube-3}, the height of primary vortex calculated by 
the current schemes agree well with the reference data \cite{Case-Daru, GKS-high-1}. 
The results also indicate the robustness of implicit methods for viscous flow with strong discontinuities.

\begin{figure}[!h]
\centering
\includegraphics[width=0.48\textwidth]{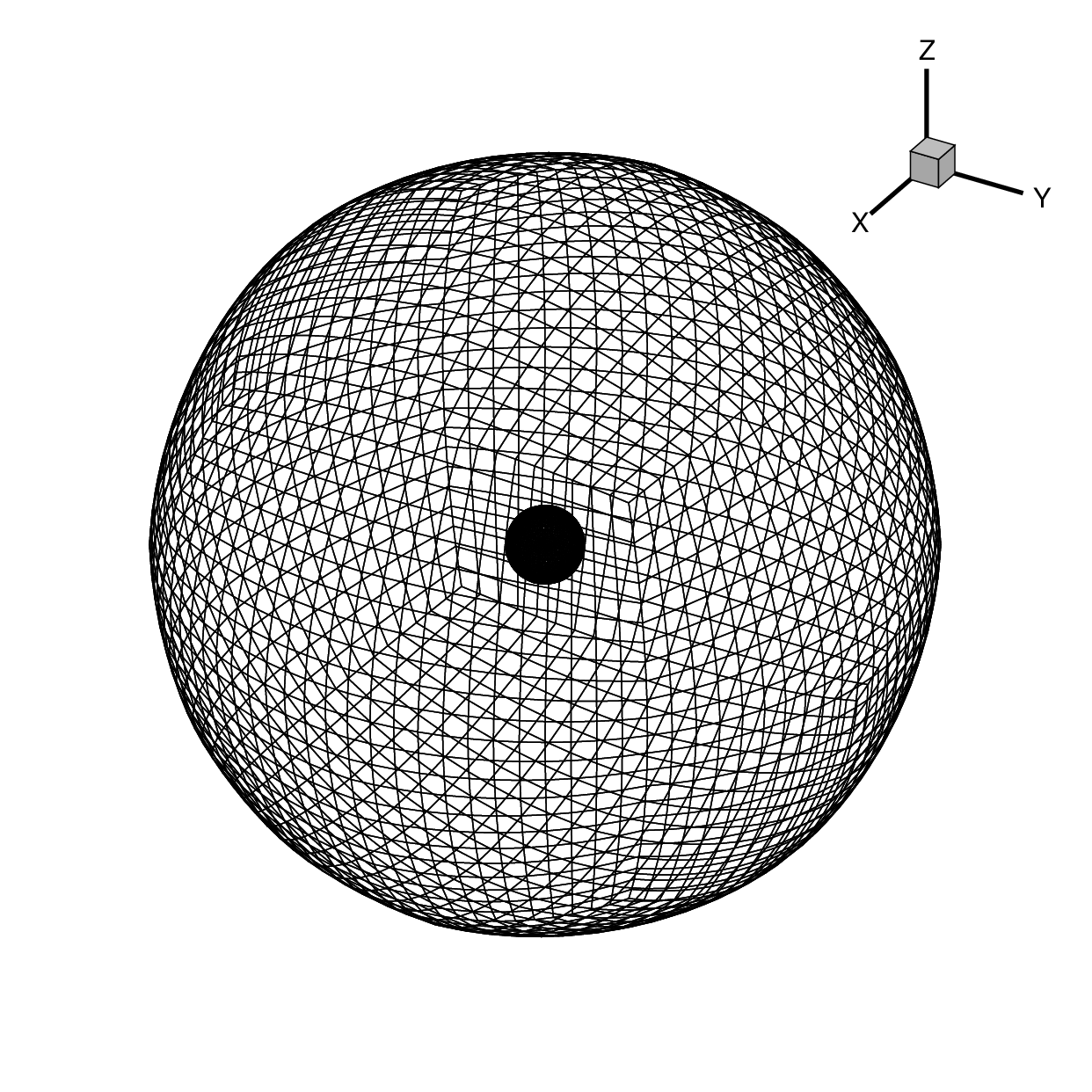}
\caption{\label{sphere-shock-mesh} Shock-sphere interaction: the computational mesh distribution.}
\end{figure}

\begin{figure}[!h]
\centering
\includegraphics[width=0.45\textwidth]{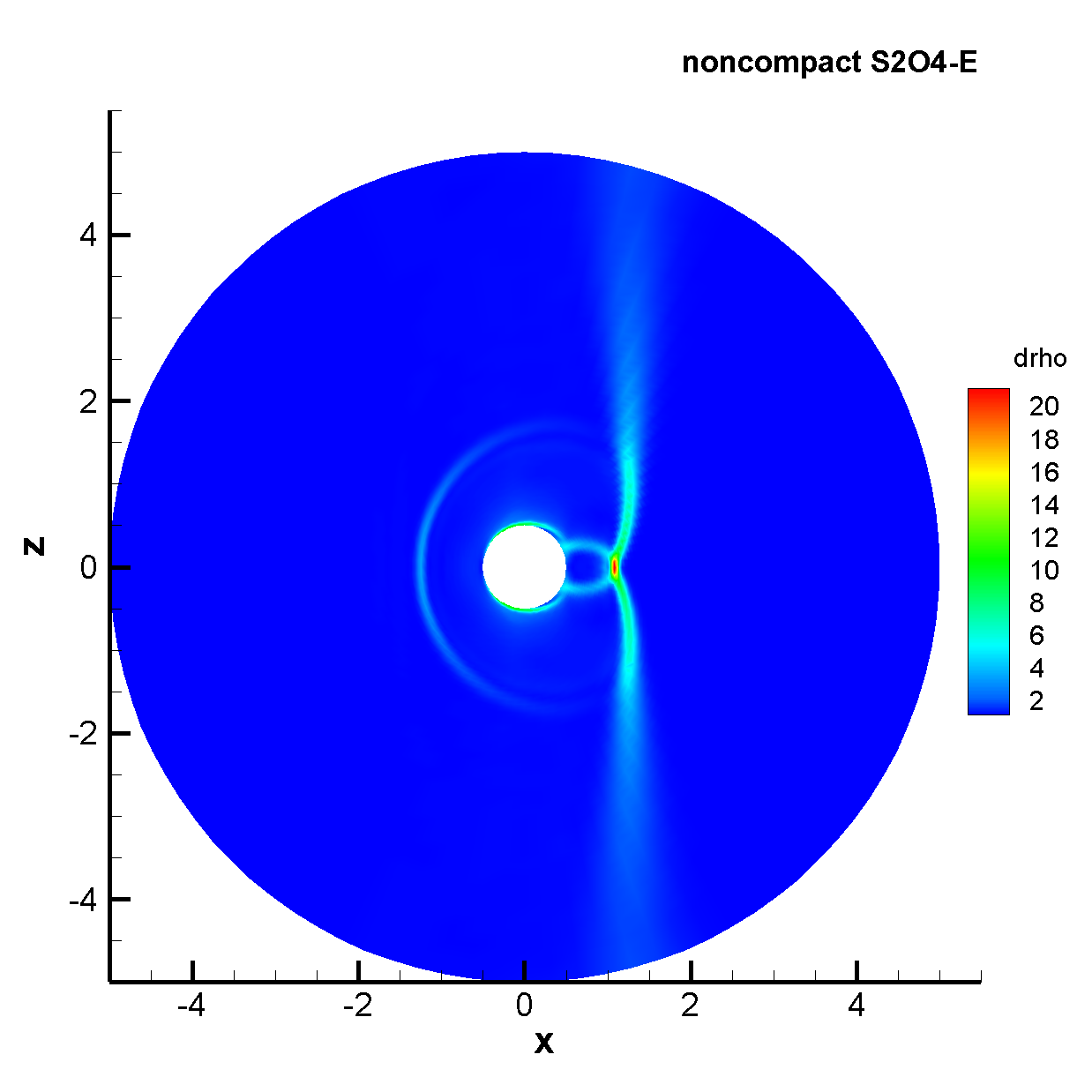}
\includegraphics[width=0.45\textwidth]{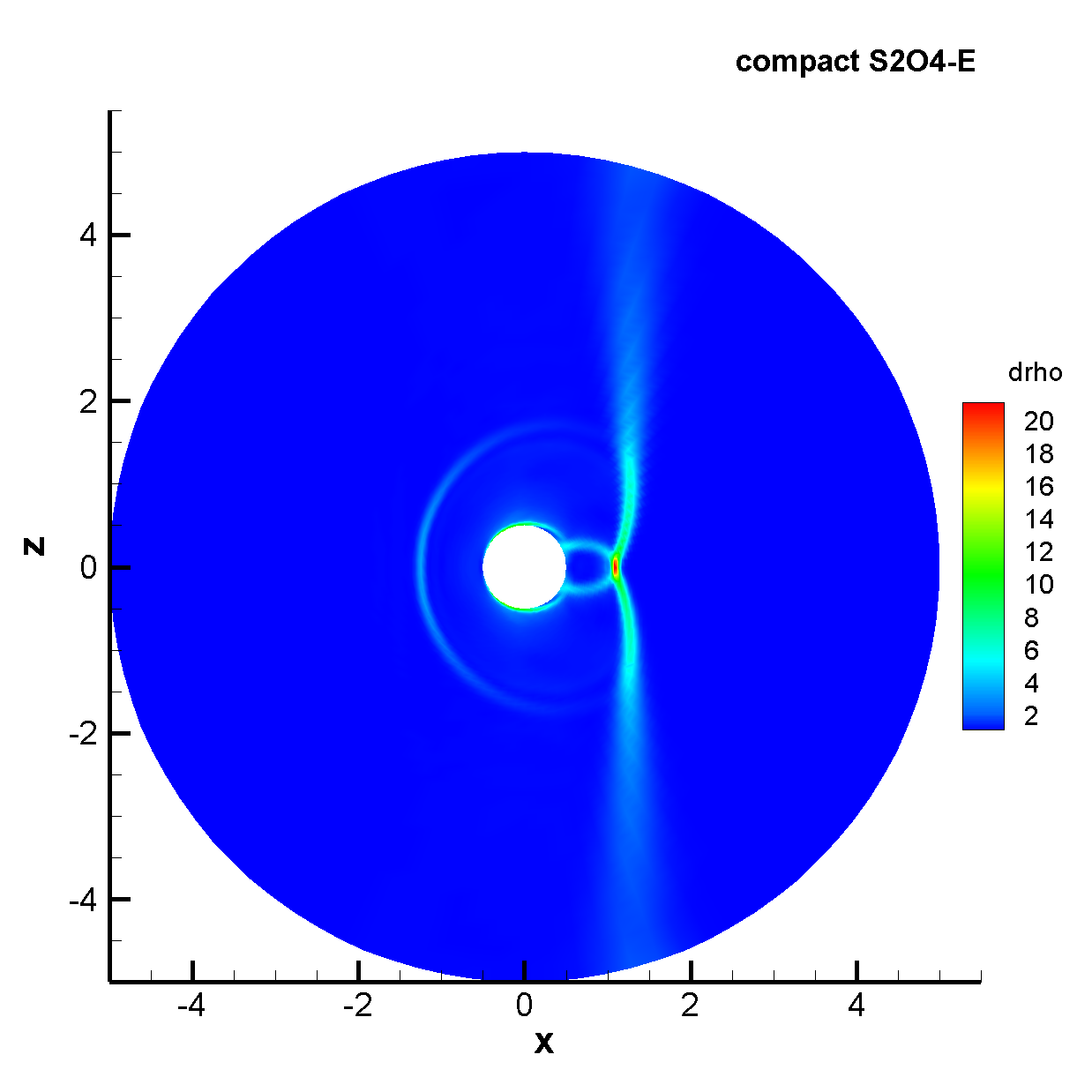}\\ 
\includegraphics[width=0.45\textwidth]{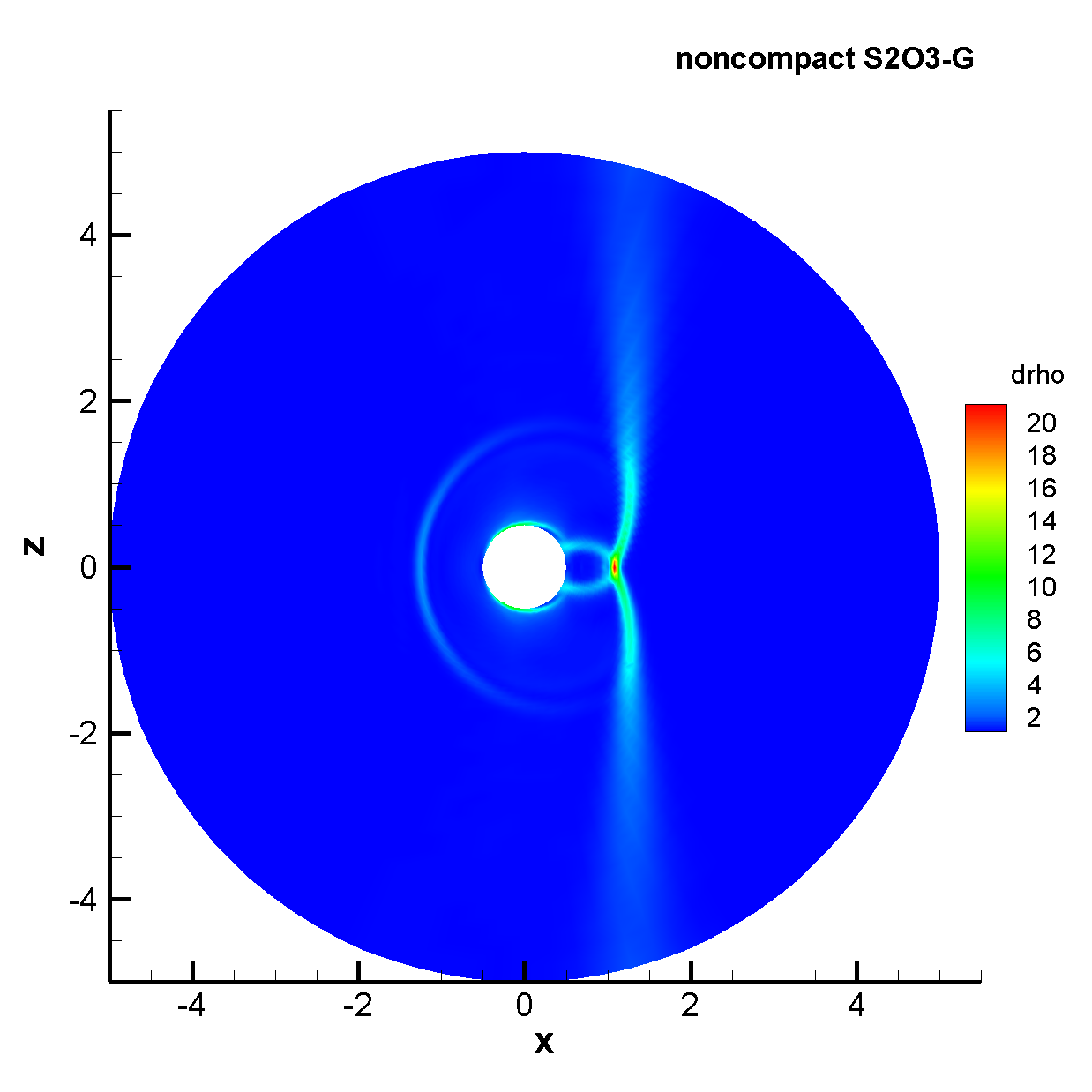} 
\includegraphics[width=0.45\textwidth]{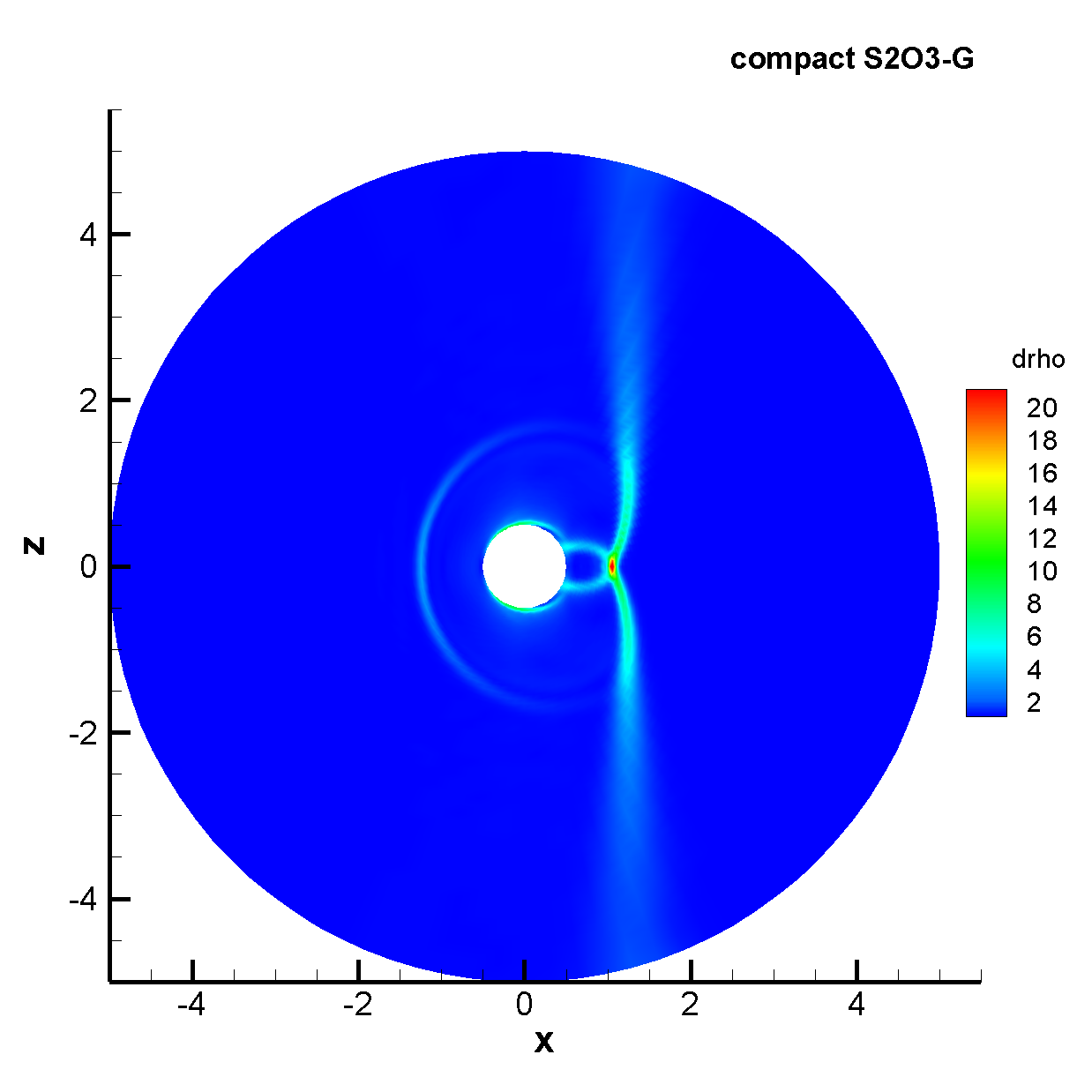} 
\caption{\label{sphere-shock-1} Shock-sphere interaction:  the density gradient distributions with non-compact and compact S2O3-G at $t = 1.5$.}
\end{figure}

\subsection{Shock-sphere interaction}
In order to explore the effect of mesh variation on the computational efficiency, a shock-sphere interaction is tested. 
This case is tested on the hexahedral mesh which contains $276480$ cells, where the variation of largest and smallest cells is about $3000$. 
The computational mesh is shown in Figure.\ref{sphere-shock-mesh}, in which the center of sphere locates at $(0,0,0)$, the radius of sphere is $R=0.5$, and 
the size of the first layer cells on the surface of sphere is $h_{\text{min}}=0.01$. The initial shock wave with Mach number $Ma=1.5$ 
and Reynolds number  $Re=1000$ locates at $X=-1.0$, and impinges on a sphere from the left. The pre-shock and post-shock states are 
given by the Rankine-Hugoniot relation. The non-slip adiabatic boundary condition is imposed for viscous flows on the surface of sphere 
and the symmetric boundary condition is used for far field.  For the implicit methods, the times of artificial iterations $k_a$ is set as 2, and 
the smallest time steps in the whole computation satisfy $\Delta t > 2\times 10^{-2}$. 
The schlieren type images of density distribution for noncompact and compact S2O3-G at $t=1.5$ are given in Figure.\ref{sphere-shock-1}.
As reference, the results with S2O4-E methods are given as well.
Figure.\ref{sphere-shock-1} shows that the numerical results computed by implicit methods match well with explicit methods.
The strong discontinuities are well resolved, which validate the robustness and accuracy of the current implicit schemes.

\begin{figure}[!h]
\centering
\includegraphics[width=0.45\textwidth]{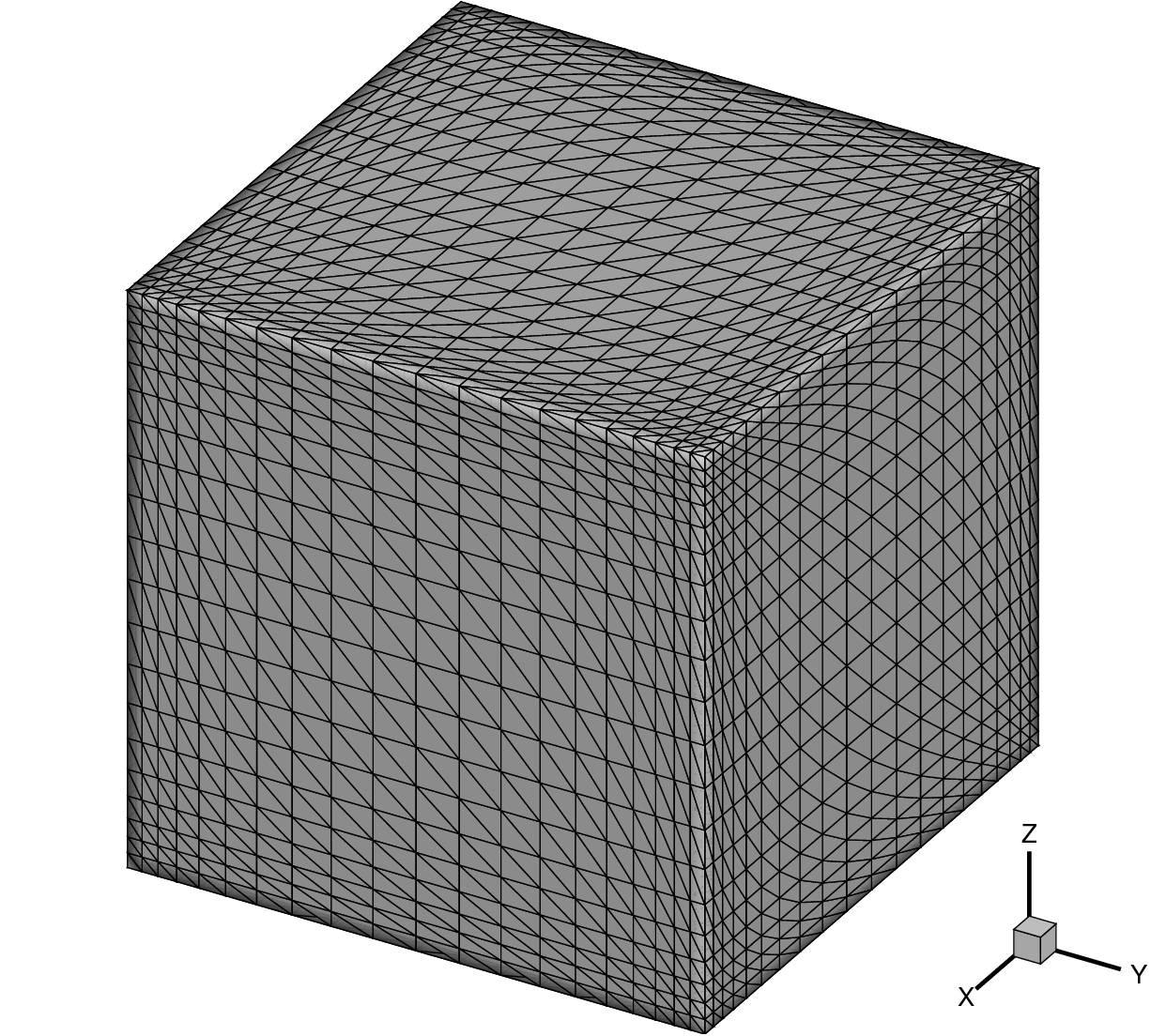}
\caption{\label{cavity-mesh} Lid-driven cavity flow: the computational mesh distribution.}
\centering
\includegraphics[width=0.425\textwidth]{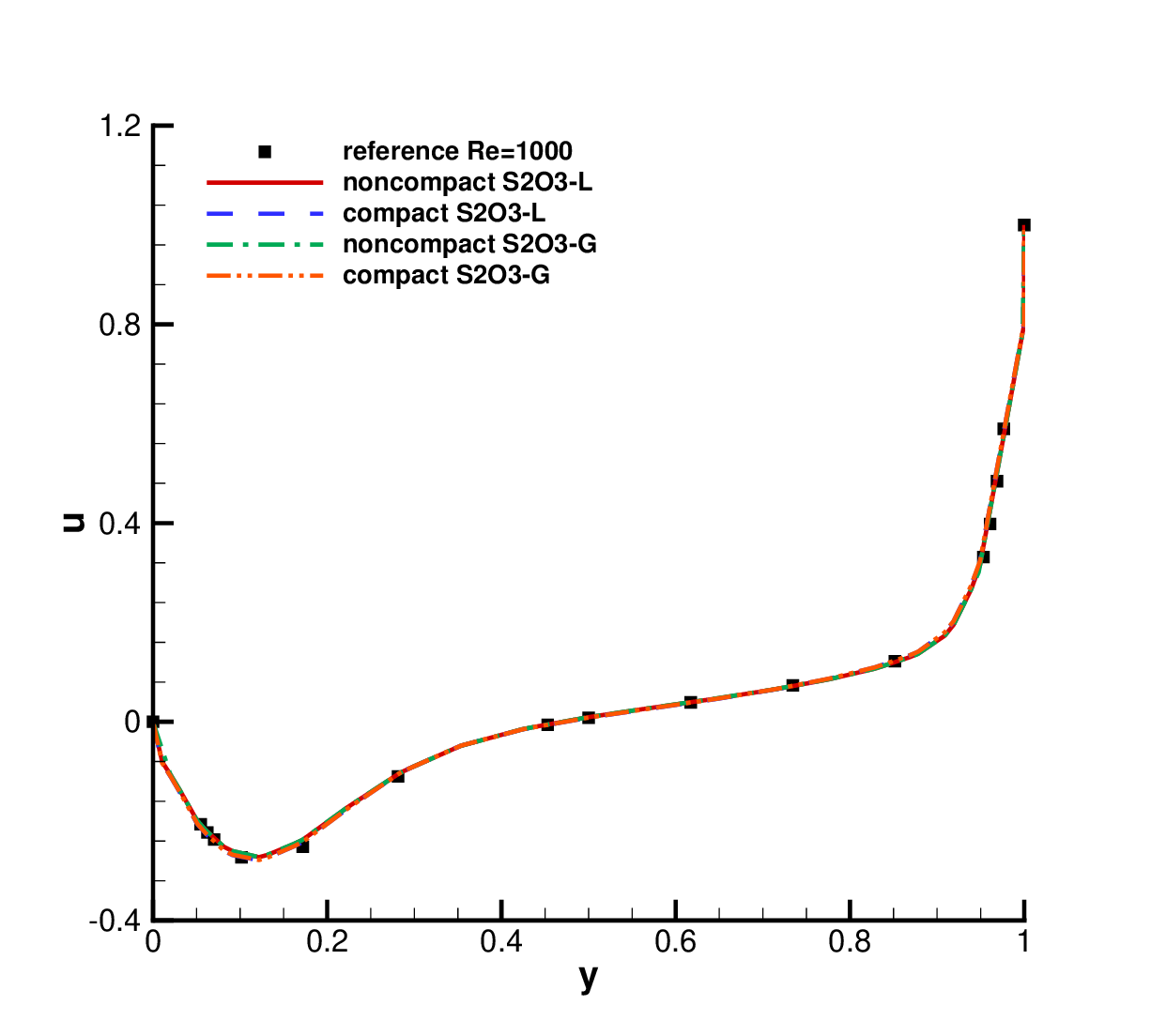}
\includegraphics[width=0.425\textwidth]{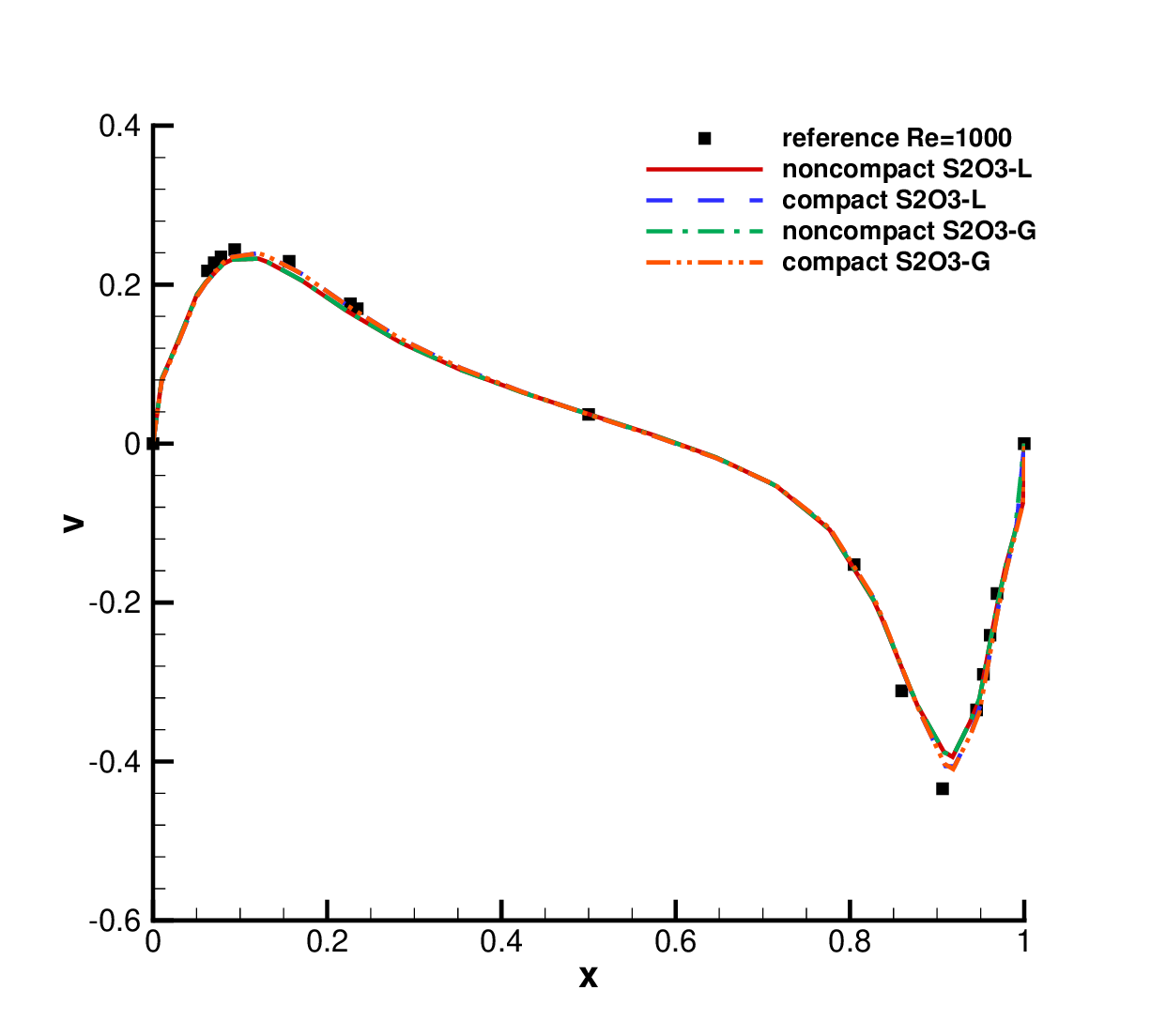}
\caption{\label{cavity-1000-1} Lid-driven cavity flow: the U-velocity profiles along the vertical centerline, V-velocity profiles along the horizontal centerline for $Re=1000$.}
\end{figure}

\subsection{Lid-driven cavity flow}
The lid-driven cavity problem is one of the most important benchmarks for numerical Navier-Stokes solvers. 
The fluid is bounded by a unit cubic $[0, 1]\times[0, 1]\times[0, 1]$ and driven by a uniform translation of the
 top boundary with $y=1$.  In this case, the flow is simulated with Mach number $Ma=0.15$ and all the boundaries 
 are isothermal and nonslip. The tetrahedral mesh with $6\times20^3$ cells is used, in which the mesh near the well 
 is refined as shown in Figure.\ref{cavity-mesh} and the size of the first layer cells is $h_{\text{min}} = 2.5\times 10^{-2}$.
Numerical simulations are conducted with the Reynolds numbers of $Re=3200$ and $1000$ by non-compact and compact S2O3-L and S2O3-G schemes. 
In this case, the times of artificial iterations $k_a$ is set as 1 for all implicit schemes.

\begin{figure}[!h]
\centering
\includegraphics[width=0.5\textwidth]{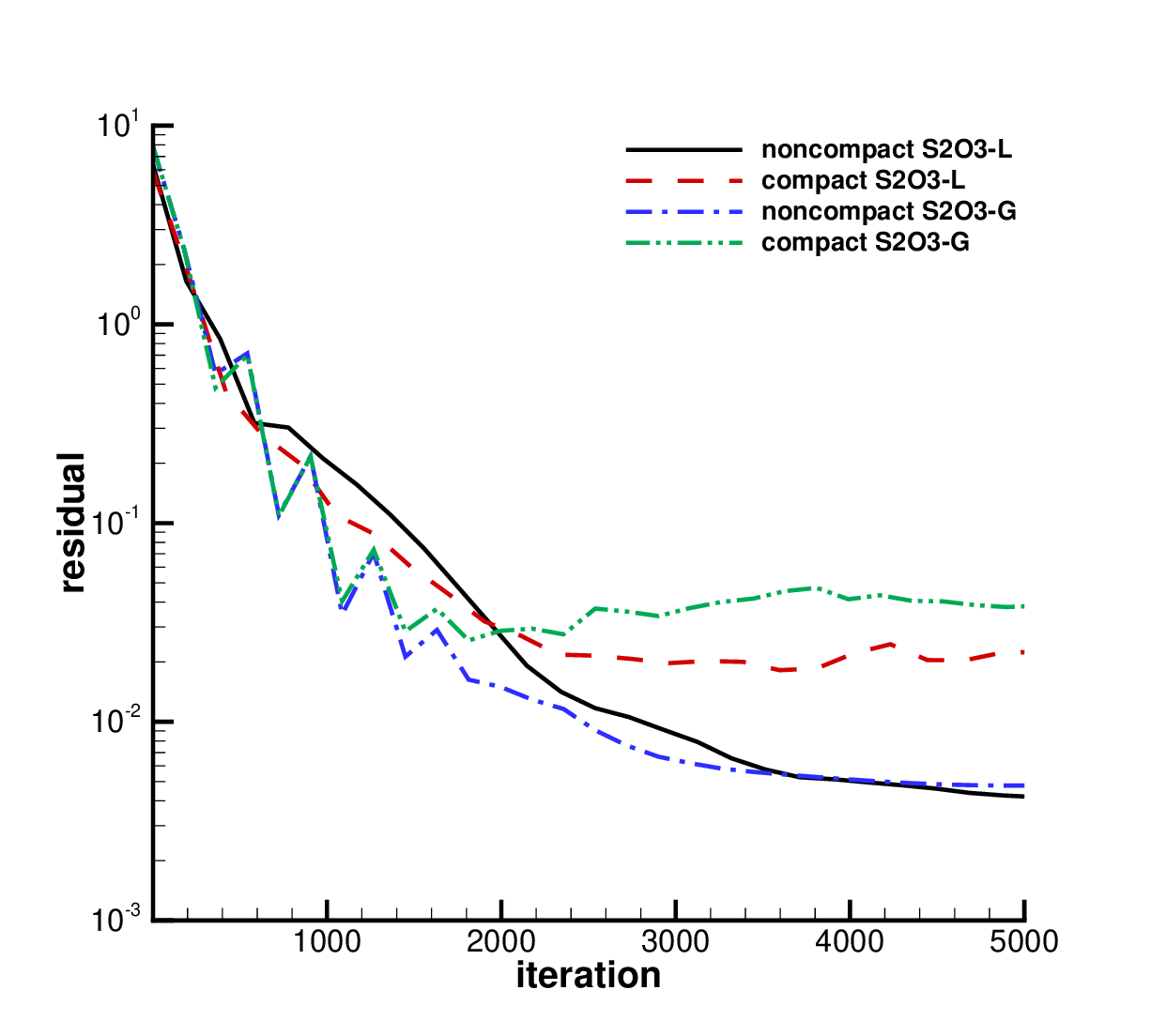}
\caption{\label{cavity-1000-2} Lid-driven cavity flow: the residual comparison for $Re=1000$.}
\end{figure}

\begin{figure}[!h]
\centering
\includegraphics[width=0.48\textwidth]{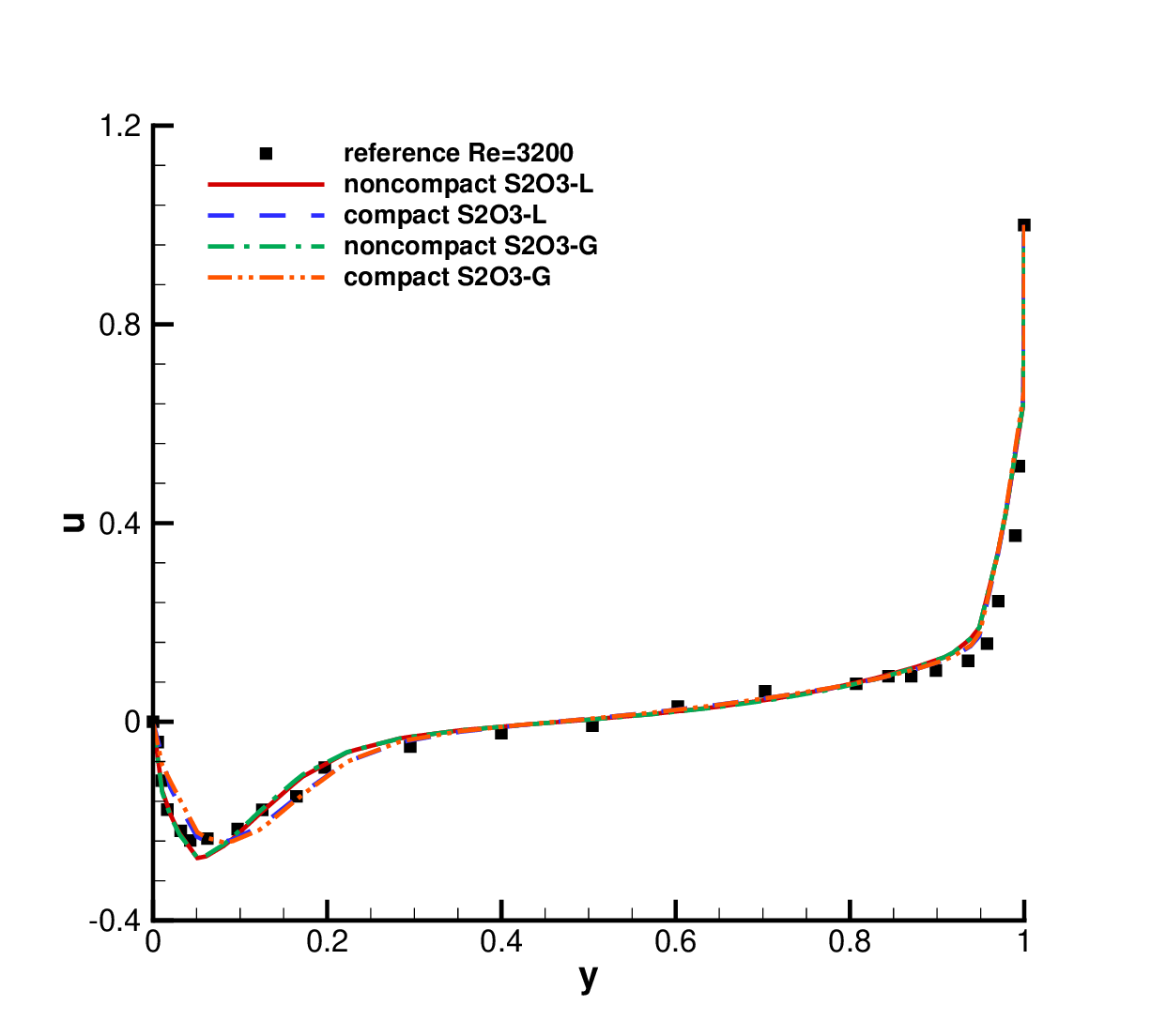}
\includegraphics[width=0.48\textwidth]{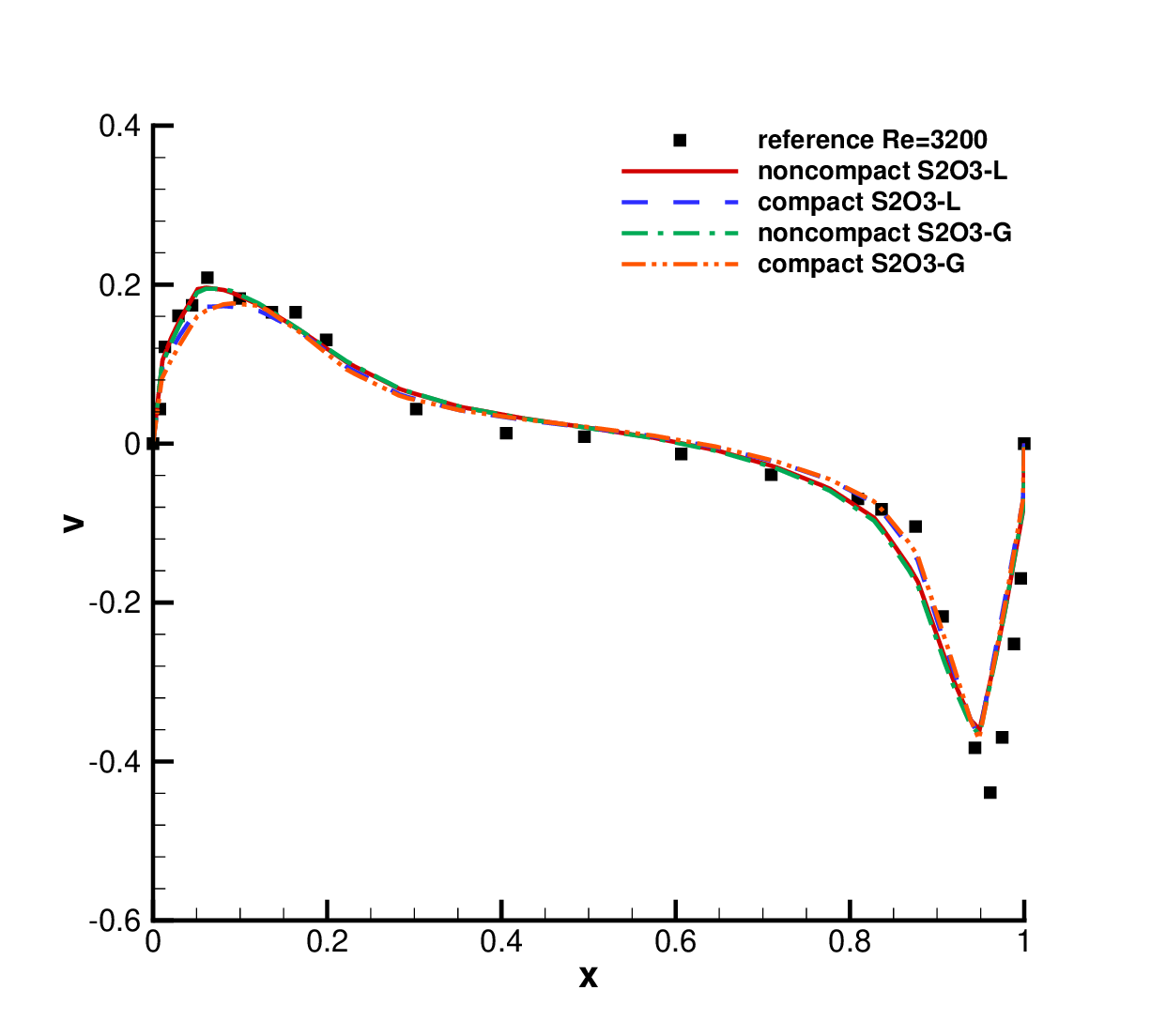}
\caption{\label{cavity-1000-3} Lid-driven cavity flow: the averaged $U$-velocity profiles along the vertical centerline, $V$-velocity profiles along the horizontal centerline for $Re=3200$.}
\end{figure}

The case with $Re=1000$ corresponds to a convergent solution. 
This case is tested to ensure that current schemes which act as unsteady 
problem solvers can simulate the steady problem correctly.  
The $U$-velocity profiles along the vertical centerline line, $V$-velocity profiles along the 
horizontal centerline in the symmetry $x-y$ plane, and the benchmark data \cite{Case-Albensoeder} are shown in Figure.\ref{cavity-1000-1}. 
The numerical results  agree well with the benchmark data, especially the compact schemes. The histories of residual 
convergence with different implicit methods are shown in Figure.\ref{cavity-1000-2}. It can be observed that 
all the implicit schemes take less than 2500 steps to converge to the steady state and the non-compact schemes have smaller converge residuals.
The case with  $Re=3200$ corresponds to an unsteady solution \cite{Case-Prasad}. 
In this case, the CFL number is set as 6.0 for all the implicit schemes and the dimension of Krylov subspace $\dim K=2$. 
The numerical results are averaged by 100 time period from $T=100$ to $200$.   
The averaged $U$-velocity profiles along the vertical centerline line, $V$-velocity profiles along the horizontal centerline in the symmetry $x-y$ plane, 
and the benchmark data \cite{Case-Prasad} are shown in Figure.\ref{cavity-1000-3}, and the numerical results also match well with the benchmark data.

\begin{table}[!h]
\begin{center}
\centering
\begin{tabular}{ccc}
\toprule
Scheme                      & CPU time &  Ratio \\
\midrule
noncompact S2O4-E &   5242s          &  ~       \\
noncompact S2O3-G &    4338s         &  0.83   \\
noncompact S2O3-L   &  5335s        &  1.02  \\
\bottomrule
\toprule
Scheme                      & CPU time &  Ratio \\
\midrule
compact S2O4-E        &  6150s        &  ~       \\
compact S2O3-G        &  4922s        &    0.80   \\
compact S2O3-L         &  6087s       &   0.99   \\
\bottomrule
\end{tabular}
\caption{\label{Efficiency-1} Efficiency comparison: the ratio of computational time for two-dimensional Riemann problem.}
\end{center}
\begin{center}
\centering
\begin{tabular}{ccc}
\toprule
Scheme                      & CPU time &  Ratio \\
\midrule
noncompact S2O4-E &   27457s          &  ~       \\
noncompact S2O3-G &    23776s         &  0.87   \\
noncompact S2O3-L   &  28974s        &  1.06  \\
\bottomrule
\toprule
Scheme                      & CPU time &  Ratio \\
\midrule
compact S2O4-E        &  33739s        &  ~       \\
compact S2O3-G        &  25318s        &    0.75   \\
compact S2O3-L         &  26354s       &   0.78   \\
\bottomrule
\end{tabular}
\caption{\label{Efficiency-2} Efficiency comparison: the ratio of computational time for viscous shock tube problem.}
\end{center}
\end{table}
\begin{table}[!h]
\begin{center}
\centering
\begin{tabular}{ccc}
\toprule
Scheme                      & CPU time &  Ratio \\
\midrule
noncompact S2O4-E &  2741s   &  ~    \\
noncompact S2O3-G &  1397s   & 0.51  \\ 
\bottomrule
\toprule
Scheme                      & CPU time &  Ratio  \\
\midrule
compact S2O4-E    & 3406s &  ~    \\
compact S2O3-G    & 1676s &  0.49     \\ 
\bottomrule
\end{tabular}
\caption{\label{Efficiency-4} Efficiency comparison: the ratio of computational time for shock-sphere interaction.}
\end{center}
\begin{center}
\centering
\begin{tabular}{ccc}
\toprule
Scheme                      & CPU time &  Ratio \\
\midrule
noncompact S2O4-E &   242s        &  ~       \\
noncompact S2O3-G &    97s         &  0.40  \\ 
\bottomrule
\toprule
Scheme                      & CPU time &  Ratio \\
\midrule
compact S2O4-E        &  335s        &  ~       \\
compact S2O3-G        &  104s       &   0.31  \\ 
\bottomrule
\end{tabular}
\caption{\label{Efficiency-3} Efficiency comparison: the ratio of computational time for lid-driven cavity flow with $Re=3200$.}
\end{center}
\end{table}

\subsection{Efficiency comparisons}
The efficiency comparison of implicit and explicit schemes is provided for both compact and non-compact reconstructions.
The CPU codes are run with AMD EPYC 9534 64 Core Processor using Intel Fortran compiler with 48 OpenMP threads, 
and the clock rate is 2450 MHz. 
\begin{enumerate}
\item For the accuracy test, the times of artificial iterations have to be set big enough for controlling the size of residual, which affects errors directly.  
For the one-dimensional Riemann problem, both implicit and explicit schemes only take few seconds to complete the computation. 
Therefore, the efficiency comparisons of these two cases are also out of consideration. 
\item For the two-dimensional Riemann problem, the efficiency comparisons are shown in 
Table.\ref{Efficiency-1}. For the viscous shock tube problem, the efficiency comparisons are 
shown in Table.\ref{Efficiency-2}. It can be observed that S2O3-G scheme has a more clear efficiency improvement 
than S2O3-L scheme for both non-compact and compact reconstructions. The reason is that the two serial sweep processes 
cells in LUSGS scheme take too much computational time with a large number of cells, and these processes are not easy to be parallelized.
For practical purpose, the S2O3-L scheme are not chosen in the following cases.
\item For the shock-sphere interaction, the efficiency comparisons at $T=[0,3]$ are shown in Table.\ref{Efficiency-4}.
Even with a shock wave, the non-compact and compact S2O3-G schemes are about two times faster than explicit schemes. 
Compared with the small acceleration in Table.\ref{Efficiency-1} and Table.\ref{Efficiency-2}, the acceleration ratio is mainly due to the large variation of mesh size.
\item
 For the lid-driven cavity flow with $Re=3200$, the efficiency comparisons within $T=[0,1]$ are shown in Table.\ref{Efficiency-3}.
In this case, the implicit schemes have obvious efficiency improvement. Especially, the compact implicit schemes are more than 
three times faster than compact S2O4-E scheme. In this case, the variation of largest and smallest cells is only about $10$. 
The main reason for such acceleration ratio is that the smooth flow only requires a single artificial iteration. Meanwhile, for 
the case with discontinuities, the increase of the  artificial iteration steps introduce extra reconstructions associating with a reduction of the computational efficiency. In order to further improve the computational efficiency,  
more effective parallel strategies need to be considered for the implicit methods.
 \end{enumerate}

\section{Conclusion}

For the simulations of unsteady compressible flow,
a new time-implicit high-order gas-kinetic scheme is developed based on the three-dimensional unstructured meshes. 
In order to enlarge the time step and keep high-order accuracy for unsteady simulations,  a two-stage third-order implicit time-accurate discretization is proposed.
In each stage, an artificial steady solution is obtained for the implicit system with the pseudo-time iteration. The LUSGS and GMRES methods are adopted to solve the implicit system in the pseudo-time iteration. 
To achieve the spatial accuracy, the third-order non-compact and compact reconstructions both are used. 
Numerical tests validate the expected order of accuracy and high efficiency in the capturing of flow structures with a large time step. Compared with the classical 
second-order Crank-Nicolson scheme, the current method achieves the higher-order temporal accuracy, and suppresses the spurious oscillations near the strong discontinuities. The numerical results also indicate that the current implicit method leads to a significant improvement of efficiency over the explicit one in the cases with a large variation of mesh size.

\section*{Acknowledgements}
The current research of L. Pan is supported by Beijing Natural
Science Foundation (1232012), National Natural Science Foundation of
China (11701038) and the Fundamental Research Funds for the Central
Universities, China.  The work of K. Xu is supported by National Key
R$\&$D Program of China (2022YFA1004500), National Natural Science
Foundation of China (12172316), and Hong Kong research grant council
(16208021,16301222).

\section*{Declaration of competing interest}
The authors declare that they have no known competing financial interests or personal relationships 
that could have appeared to influence the work reported in this paper.

\section*{Data availability}
The data that support the findings of this study are available from
the corresponding author upon reasonable request.

\end{document}